\documentclass[11pt, final]{article}
\usepackage{a4}
\usepackage{amsmath}%
\usepackage{amstext}%
\usepackage{amssymb}%
\usepackage{showkeys}%
\usepackage{epsfig}%
\usepackage{cite}
\newtheorem{theorem}{Theorem}

\newtheorem{example}[theorem]{Example}

\newtheorem{lemma}[theorem]{Lemma}

\newtheorem{proposition}[theorem]{Proposition}
\newtheorem{remark}{Remark}

\newcounter{unnumber}

\newcommand{\R}{\mathbb{R}}%
\newcommand{\N}{\mathbb{N}}%
\DeclareMathOperator*\inte{int}%
\DeclareMathOperator*\ri{ri}%
\DeclareMathOperator*\cl{cl}%
\DeclareMathOperator*\cp{\cup}%
\DeclareMathOperator*\epi{epi}%
\DeclareMathOperator*\dom{dom}%
\DeclareMathOperator*\aff{aff}
\DeclareMathOperator*\qi{qi}%
\DeclareMathOperator*\qri{qri}%
\DeclareMathOperator*\sqri{sqri}%
\DeclareMathOperator*\epihat{\widehat{\epi}}%
\DeclareMathOperator*\coneco{coneco}%
\DeclareMathOperator*\cone{cone}%
\DeclareMathOperator*\core{core}%
\DeclareMathOperator*\icr{icr}%
\DeclareMathOperator*\co{co}%
\DeclareMathOperator*\B{\overline{\R}}%
\DeclareMathOperator*\id{id}%
\DeclareMathOperator*\pr{pr}%
\DeclareMathOperator*\lin{lin}%
\DeclareMathOperator*\Ker{ker}%

\title{Regularity conditions via generalized interiority notions in convex optimization: new achievements and their relation to some classical statements}

 \author{Radu Ioan Bo\c{t} \thanks
 {Faculty of Mathematics, Chemnitz University of Technology,
D-09107 Chemnitz, Germany, e-mail:
 radu.bot@mathematik.tu-chemnitz.de. Research partially supported by DFG (German Research Foundation), project WA 922/1-3.} \and Ern\"{o} Robert Csetnek
 \thanks {Faculty of Mathematics, Chemnitz University of Technology,
D-09107 Chemnitz, Germany, e-mail:
 robert.csetnek@mathematik.tu-chemnitz.de}}

\begin{document}
\maketitle

\noindent \textbf{Abstract.} For the existence of strong duality in convex optimization regularity conditions play an indisputable role. We mainly deal in this paper
with regularity conditions formulated by means of different generalizations of the notion of interior of a set. The primal-dual pair we investigate is a general one expressed in the language
of a perturbation function and by employing its Fenchel-Moreau conjugate. After providing an overview on the generalized interior-point conditions that exist in the literature we introduce several new ones formulated by means of the quasi interior and quasi-relative interior. We underline the advantages of the new conditions vis-\' a-vis the classical ones and illustrate our investigations by numerous examples. We close the paper by particularizing the general approach to the classical Fenchel and Lagrange duality concepts.\\



\section{Introduction}

One of the most important issues that occur in the investigations made in connection to the duality theory in convex optimization is the formulation and the verification of so-called \emph{regularity conditions}. Two main classes of such conditions exist in the literature, on the one hand, the meanwhile classical \emph{generalized interior-point regularity conditions} (see \cite{Gowda-Teboulle, Rock-conj-dual, Zal-art, Zal-carte}) and, on the other hand, the recently introduced \emph{closedness-type regularity conditions} (a comprehensive reference to the latter is \cite{b-hab}). In this paper we mainly deal with regularity conditions belonging to the first class, the closedness-type ones being only tangentially addressed. To the class of generalized interior-point regularity conditions belong the ones that assume in a more or less direct manner continuity for the functions involved. On the other hand, one also has here assumptions that ask for completeness for the underlying spaces, lower semicontinuity (or generalizations of this notion) for the functions involved as well as some conditions expressed by means of the \emph{interior} but also of some of its generalizations. In this sense we mention here the \emph{algebraic interior}, the \emph{relative algebraic interior} and the \emph{strong quasi-relative interior}.

One of the aims of this paper is to provide new generalized interior-point regularity conditions this time expressed via the \emph{quasi interior} and the \emph{quasi-relative interior}. The latter
has been introduced by Borwein and Lewis in \cite{Borwein-Lewis} and constitutes a generalization of the interiority notions mentioned above. When dealing for instance with the \emph{positive
cones} of the spaces $\ell^p$ or $L^p(T,\mu)$ (here $(T,\mu)$ is a $\sigma$-finite measure space) when $p \in [1, \infty)$, we have that their quasi-relative is nonempty, unlike it happens for their interior, algebraic interior, strong quasi-relative interior and relative algebraic interior. This fact guarantees a wider applicability of the regularity conditions considered here vis-\' a-vis the classical ones, another supporting argument for the new ones being given by the fact that no completeness for the underlying spaces neither topological requirements for the functions involved are assumed. We also address the relations between the new conditions and the classical ones and illustrate our investigations by numerous examples.

In the section that follows this introductory one we establish the setting in which we work and we recall some notions and results from the convex analysis. In the third section we give first an overview on the generalized interiority notions introduced in the literature in connection with a convex set, followed by the definitions of the quasi interior and the quasi-relative interior. For the latter we furnish dual characterizations and list some essential properties. A scheme with the existing relations between the addressed generalized interiority notions is also provided, the section being closed by a separation result formulated by employing the quasi-relative interior. In section 4 we deal with a general \emph{conjugate duality} scheme for convex optimization problems having as starting point the meanwhile classical \emph{perturbation approach} (see for instance \cite{Rock-conj-dual, EkTem, Zal-carte}). Some of the generalized interior-point regularity conditions that have been introduced in the literature in this general context are recalled, followed by the formulation of two new conditions expressed via the \emph{quasi interior} and the \emph{quasi-relative interior}. After proving that the latter are sufficient for having strong duality we also succeed to relate them to the classical ones. The general duality scheme is particularized in the next section to the problem having the sum of two convex functions as objective function and its \textit{Fenchel dual problem}. Here we also formulate a supplementary regularity condition strongly connected to the two new ones.
Further, we illustrate via several examples the applicability of the new conditions. Moreover, we show that they cannot be compared for the closedness-type conditions considered in the literature for the same primal-dual pair and give in this was a negative answer to an open problem stated in \cite{li-fang-lopez}. We close the paper by dealing with the optimization problem with geometric and cone constraints and its \textit{Lagrange dual problem}, again seen as a particular instance of the general primal-dual pair from section 4. We do not only particularizing the two general regularity conditions, but also formulate a further one strongly connected to them. Also here we illustrate the relations between all these conditions as well as the applicability of the new ones by several examples.

\section{Preliminaries}

Consider $X$ a separated locally convex space and $X^*$ its topological dual space. We denote by $w(X^*,X)$ the weak$^*$ topology on $X^*$ induced by $X$. For a nonempty set $U\subseteq X$, we denote by $\co(U), \cone(U), \coneco(U), \aff(U), \lin(U), \inte(U), \cl(U)$, its \emph{convex hull, conic  hull, convex conic hull,  affine hull, linear hull, interior} and \emph{closure}, respectively. We have $\cone (U)=\cup_{t\geq 0}t U$ and if $0\in U$ then obviously $\cone (U)= \cup_{t > 0}t U$. In case $U$ is a linear subspace of $X$ we denote by $U^\bot$ the \emph{annihilator} of $U$. The following property will be used several times throughout this paper: if $U$ is convex then
\begin{equation}\coneco(U\cup\{0\})=\cone(U).\label{coneco}\end{equation} If $U\subseteq \R^n$ ($n\in\N$) we denote by $\ri(U)$ the \emph{relative interior} of $U$, that is the interior of $U$ with respect to its affine hull. We denote by
$\langle x^*,x\rangle$ the value of the linear continuous
functional $x^*\in X^*$ at $x\in X$ and by $\Ker x^*$ the \emph{kernel} of $x^*$.
Let us consider $V\subseteq X$ another nonempty set.
The \emph{projection operator} $\pr_U : U\times V\rightarrow U$ is defined as $\pr_U(u,v)=u$ for all $(u,v)\in U\times V$, while the \emph{indicator function} of $U$,
$\delta_U:X\rightarrow\B$, is defined as
$$\delta_U(x)=\left\{
\begin{array}{ll}
0, & \mbox {if } x\in U,\\
+\infty, & \mbox{otherwise},
\end{array}\right.$$ where $\B=\R\cup\{\pm\infty\}$ is the extended real line. We say that the function $f:X\rightarrow \B$ is \emph{convex} if \begin{equation}\forall x,y\in X, \ \forall t\in[0,1]: f(t x+(1-t)y)\leq t f(x)+(1-t)f(y),\end{equation} with the conventions $(+\infty)+(-\infty)=+\infty,$ $0\cdot(+\infty)=+\infty$ and $0\cdot(-\infty)=0$. We denote by $\dom f=\{x\in X:f(x)<+\infty\}$
the \emph{domain} of $f$ and by $\epi f=\{(x,r)\in X\times\R:f(x)\leq
r\}$ its \emph{epigraph}. Moreover, we denote by $\epihat(f)=\{(x,r)\in X\times\R:(x,-r)\in\epi f\}$, the
symmetric of $\epi f$ with respect to the $x$-axis. For a given
real number $\alpha$, $f-\alpha:X\rightarrow\B$ is, as usual, the
function defined by $(f-\alpha)(x)=f(x)-\alpha$ for all $x\in X$. We call $f$ \emph{proper} if
$\dom f\neq\emptyset$ and $f(x)>-\infty$ for all $x\in X$. The \emph{normal cone} of $U$ at $x\in U$ is $N_U(x)=\{x^*\in X^*:\langle x^*,y-x\rangle\leq 0 \ \forall y\in U\}$.

The \emph{Fenchel-Moreau conjugate} of $f$ is the function
$f^*:X^*\rightarrow\B$ defined by
$$f^*(x^*)=\sup\limits_{x\in X}\{\langle
x^*,x\rangle-f(x)\} \ \forall x^*\in X^*.$$ We have the so called \emph{Young-Fenchel
inequality}
$$f^*(x^*)+f(x)\geq\langle x^*,x\rangle \ \forall x\in X \ \forall
x^*\in X^*.$$ Having $f,g:X\rightarrow\B$ two functions we denote by
$f\Box g:X\rightarrow\B$, $f\Box g(x)=\inf_{u\in X}\{f(u)+g(x-u)\}$ for
all $x\in X$, their \emph{infimal convolution}. We say that the infimal convolution is \emph{exact
at $x\in X$} if the infimum in its definition is attained.
Moreover, $f\Box g$ is said to be \emph{exact} if it is exact at
every $x\in X$.

Consider $Y$ another separated locally convex space. For a vector function $h:X\rightarrow Y$ we denote by $h(U)=\{h(u):u\in U\}$ the \emph{image} of the set $U\subseteq X$ through $h$, while $h^{-1}(D)=\{x\in X:h(x)\in D\}$ is the \emph{counter image} of the set $D\subseteq Y$ through $h$. Given a linear continuous mapping $A:X\rightarrow Y$, its \emph{adjoint operator} $A^*:Y^*\rightarrow X^*$ is defined by
$\langle A^*y^*,x\rangle=\langle y^*,Ax\rangle$ for all $y^*\in
Y^*$ and $x\in X$. For a nonempty \emph{convex cone} $C\subseteq Y$ (that is $\cone (C)\subseteq C$ and $C+C\subseteq C$) we denote by
$C^*=\{y^*\in Y^*:\langle y^*,y\rangle\geq 0 \ \forall y\in C\}$
its \emph{positive dual cone}. Further, we denote by $\leq_C$ the partial ordering induced by $C$ on $Y$, defined as $y_1\leq_C
y_2\Leftrightarrow y_2-y_1\in C$ for $y_1,y_2\in Y$. To $Y$ we attach an
abstract maximal element with respect to $\leq_C$, denoted by
$\infty_C$ and let $Y^\bullet:=Y\cup\{\infty_C\}$. Then for every
$y\in Y$ one has $y\leq_C\infty_C$, while on $Y^\bullet$
the following operations are considered: $y+\infty_C=\infty_C+y=\infty_C$ and
$t\infty_C=\infty_C$ for all $y\in Y$ and all $t\geq 0$. Moreover, if $\lambda\in C^*$ let $\langle\lambda,\infty_C\rangle:=+\infty$.

Some of the above notions given for functions with extended real
values can be formulated also for function having their ranges in
infinite-dimensional spaces.

For a function $h:X\rightarrow Y^\bullet$ we denote by
$\dom h=\{x\in X:h(x)\in Y\}$ its \emph{domain} and by
$\epi_C h=\{(x,y)\in X\times Y:h(x)\leq_C y\}$ its
\emph{C-epigraph}. We say that $h$ is \emph{proper} if its domain
is a nonempty set. The function $h$ is said to be
\emph{$C$-convex} if
$$\forall x,y\in X, \ \forall t\in[0,1]: h(t x+(1-t) y)\leq_C th(x)+(1-t)h(y).$$ One can prove that $h$ is $C$-onvex if and only if $\epi_C h$ is a convex subset of $X\times Y$. Further, for an arbitrary $y^*\in C^*$ we define the function $(y^* h): X \rightarrow\B$, by $(y^* h)(x)=
\langle y^*,h(x)\rangle$ for all $x\in X$. The function $h$ is said to be
\emph{$C$-epi closed} if $\epi_C h$ is a closed subset of
$X\times Y$ (cf. \cite{Luc}). Another generalization of the
lower semicontinuity for functions taking values in infinite-dimensional spaces is the so-called \emph{star C-lower
semicontinuity}, namely $h$ is called \emph{star C-lower
semicontinuous} at $x\in X$ if for all $y^* \in C^*$ the
function $(y^*h)$ is lower semicontinuous at $x$. The function
$h$ is said to be \emph{star C-lower semicontinuous} if it is star
$C$-lower semicontinuous at every $x\in X$. This notion was considered first in \cite{jeyak-prep}. Let us mention that there exists in the literature another notion of lower semicontinuity, called \emph{$C$-lower semicontinuity}, which has been introduced by Penot and Th\'{e}ra in \cite{PenotThera} and then
refined in \cite{Combari-Laghdir-Thibault}. One can show that
$C$-lower semicontinuity implies star $C$-lower semicontinuity,
which yields $C$-epi-closedness (see \cite{Luc}), while the
opposite assertions are not valid in general. An example of a $C$-convex function which is $C$-epi-closed, but not star $C$-lower semicontinuous is given in \cite[Example 1]{BGW}. For more on the lower semicontinuity for functions with values in topological vector spaces we refer the reader to \cite{ManMThera, Combari-Laghdir-Thibault, Luc, PenotThera, Tanaka}. It is known that in case $Y=\R$ and $C=\R_+$, all the lower semicontinuity notions mentioned above coincide, becoming the classical lower semicontinuity for  functions with extended real values.

\section{Generalized interiority notions}

In this section we recall the most important generalized interiority notions one can find in the literature. Consider $X$ a separated locally convex space and $U\subseteq X$ a nonempty convex set. We have:

\begin{enumerate}

\item[$\bullet$] $\core(U):=\{x\in U:\cone(U-x)=X\}$, the \emph{algebraic interior} (the \emph{core}) of $U$ (cf. \cite{Rock-conj-dual, Zal-carte});

\item[$\bullet$]$\icr(U):=\{x\in U:\cone(U-x)\mbox{ is a linear
subspace of }X\}$, the \emph{relative algebraic interior}
   (\emph{intrinsic core}) of $U$ (cf. \cite{Borwein-Goebel, holmes, Zal-carte});

\item[$\bullet$] $\sqri(U):=\{x\in U:\cone(U-x)\mbox{ is a closed
linear subspace of }X\}$ the \emph{strong quasi-relative interior}
     (\emph{intrinsic relative algebraic interior}) of $U$ (cf. \cite{Borwein-Jeyak-Lewis-Wolk, Zal-carte, Jeyak-Wolk}).

\end{enumerate}

We mention the following characterization of the strong quasi-relative interior (cf. \cite{Gowda-Teboulle, Zal-carte}):
\begin{equation}\label{charicr}
x\in\sqri(U) \Leftrightarrow x\in\icr(U) \ \mbox{and} \ \aff(U) \  \mbox{is a closed linear subspace of} \ X.
\end{equation}

The \emph{quasi-relative interior} of $U$ is the set (cf.
\cite{Borwein-Lewis}) $$\qri(U)=\{x\in
U:\cl\big(\cone(U-x)\big)\mbox{ is a linear subspace of }X\}.$$

A useful characterization of the quasi-relative interior of a convex set by means of the normal cone follows.

\begin{proposition}\label{charact-qri} (cf. \cite{Borwein-Lewis}) {\it Let $U$ be
a nonempty convex subset of $X$ and $x\in U$. Then $x\in\qri(U)$
if and only if $N_U(x)$ is a linear subspace of $X^*$.}\end{proposition}

Next we consider another generalized interior-notion introduced in connection with a convex
set, which is close to the quasi-relative interior. The
\emph{quasi interior} of $U$ is the set $$\qi(U)=\{x\in
U:\cl\big(\cone(U-x)\big)=X\}.$$ It can be characterized as follows.

\begin{proposition}\label{charact-qi} (cf. \cite[Proposition 2.4]{qri-siam}) {\it Let $U$ be a nonempty convex subset of $X$ and $x\in U$. Then $x\in\qi(U)$ if and only if $N_U(x)=\{0\}$.}\end{proposition}

\begin{remark}\label{normed-lcs}\rm The above characterization of the quasi interior of a convex set was given in \cite{Daniele-Giuffre-Maugeri}, where the authors supposed that $X$ is a reflexive Banach space. It is proved in \cite[Proposition 2.4]{qri-siam} that this property holds in a more general context, namely in separated locally convex spaces.\end{remark}

The following scheme furnishes the relations between the different generalized interiority notions considered above
\begin{equation}\label{int-not-1}
\inte(U)\subseteq\core(U)\subseteq
\begin{array}{c}
\sqri(U)\subseteq \icr(U)\\
~\\
\qi(U)
\end{array}
\subseteq\qri(U)\subseteq U,
\end{equation}
all the inclusions being in general strict. As one can also deduce from some of the examples which follows in this paper in general between $\sqri(U)$ and $\icr(U)$, on the one hand, and $\qi(U)$, on the other hand, no relation of inclusion can be provided. In case $\inte(U)\neq\emptyset$ all the
generalized interior-notions considered in \eqref{int-not-1} collapse into
$\inte(U)$ (cf. \cite[Corollary 2.14]{Borwein-Lewis}).

It follows from the definition of the quasi-relative interior that
$\qri(\{x\})=\{x\}$ for all $x\in X$. Moreover, if $\qi(U)\neq\emptyset$, then $\qi
(U)=\qri(U)$. Although this property is given in
\cite{Limber-Goodrich} in the case of real normed spaces, it holds
also in separated locally convex spaces, as it easily follows from
the properties given above. For $U,V$ two convex subsets of $X$
such that $U\subseteq V$, we have $\qi(U)\subseteq \qi(V)$, a
property which is no longer true for the quasi-relative interior (however it holds in case $\aff(U)=\aff(V)$, see \cite[Proposition 1.12]{Cammaroto2}). If $X$ if finite-dimensional then
$\qri(U)=\sqri(U)=\icr (U)=\ri(U)$ (cf. \cite{Borwein-Lewis, Gowda-Teboulle}) and $\core(U)=\qi(U)=\inte(U)$ (cf. \cite{Limber-Goodrich, Rock-conj-dual}). We refer the reader to
\cite{Borwein-Goebel, Borwein-Lewis, Gowda-Teboulle, holmes,
Limber-Goodrich, Rock-conj-dual, Tan-Kuroiw, Zal-carte} and to the references therein for more
properties and examples regarding the above considered generalized
interiority notions.

\begin{example}\label{qri-lp+} {\rm Take an arbitrary $p\in[1,+\infty)$ and consider
the real Banach space $\ell^p=\ell^p(\N)$ of real sequences $(x_n)_{n\in\N}$ such
that $\sum\limits_{n=1}^{\infty}|x_n|^p<+\infty$, equipped with
the norm $\|\cdot\|:\ell^p\rightarrow\R$,
$\|x\|=\Big(\sum\limits_{n=1}^{\infty}|x_n|^p\Big)^{1/p}$ for all
$x=(x_n)_{n\in\N}\in \ell^p$. Then (cf. \cite{Borwein-Lewis})
$$\qri(\ell^p_+)=\{(x_n)_{n\in\N}\in \ell^p:x_n> 0 \ \forall n\in\N\},$$
where $\ell^p_+=\{(x_n)_{n\in\N}\in \ell^p:x_n\geq 0 \ \forall n\in\N\}$
is the positive cone of $\ell^p$. Moreover, one can prove that
$$\inte(\ell^p_+)=\core(\ell^p_+)=\sqri(\ell^p_+)=\icr(\ell^p_+)=\emptyset.$$}\end{example}

In a separable Banach spaces every nonempty closed convex set has a nonempty quasi-relative interior (cf. \cite[Theorem 2.19]{Borwein-Lewis}, see also \cite[Theorem 2.8]{Borwein-Goebel} and \cite[Proposition 1.2.9]{Zal-carte}) and every nonempty convex set which is not contained in a hyperplane possesses a nonempty quasi interior (cf. \cite{Limber-Goodrich}). This result may fail if we renounce to the condition $X$ is separable, as the following example shows.

\begin{example}\label{qri-vid} \rm For $p\in[1,+\infty)$ consider the real Banach space \begin{center}
\begin{tabular}{r|l}
$\ell^p(\R)=\{s:\R\rightarrow \R$ \ & \ $\sum\limits_{r\in\R} |s(r)|^p<\infty\}$,
\end{tabular}
\end{center} equipped with the norm $\|\cdot\|:\ell^p(\R)\rightarrow \R$, $\|s\|=\Big(\sum_{r\in\R}|s(r)|^p\Big)^{1/p}$ for all $s\in \ell^p(\R)$, where $$\sum_{r\in\R}|s(r)|^p=\sup\limits_{F\subseteq \R, F \mbox{\tiny finite}}\sum_{r\in F}|s(r)|^p.$$ Considering the positive cone $\ell^p_+(\R)=\{s\in\ell^p(\R):s(r)\geq 0 \ \forall r\in \R\}$, we have (cf. \cite[Example 3.11(iii)]{Borwein-Lewis}, see also \cite[Remark 2.20]{bor-luc-boris}) that $\qri\big(\ell^p_+(\R)\big)=\emptyset$.\end{example}

Some useful properties of the quasi-relative interior are listed below.
For the proof of $(i)-(ii)$ we refer to \cite{Borwein-Goebel, Borwein-Lewis}, while property $(iii)$ was proved in \cite[Proposition 2.5]{qri-siam} (see also \cite[Proposition 2.3]{qri-jota}).

\begin{proposition}\label{prop-qri} {\it Consider $U$ a nonempty convex
subset of $X$. Then:
\begin{enumerate} \item[(i)] $t\qri(U)+(1-t)U\subseteq\qri(U)$
$\forall t\in(0,1]$; hence $\qri(U)$ is a convex set.\end{enumerate} If, additionally, $\qri(U)\neq\emptyset$ then:
\begin{enumerate}\item[(ii)] $\cl\big(\qri(U)\big)=\cl(U)$; \item[(iii)] $\cl\Big(\cone\big(\qri(U)\big)\Big)=\cl\big(\cone (U)\big)$.\end{enumerate}}\end{proposition}

\noindent The first part of the next lemma was proved in \cite[Lemma 2.6]{qri-siam} (see also \cite[Lemma 2.1]{qri-jota}).

\begin{lemma}\label{alte-prop-qri} {\it Let $U$ and $V$ be nonempty convex subsets
of $X$ and $x\in X$. Then:\begin{enumerate} \item[(i)] if
$\qri(U)\cap V\neq\emptyset$ and $0\in\qi(U-U)$, then
$0\in\qi(U-V)$; \item[(ii)] $x\in\qi(U)$ if and only if
$x\in\qri(U)$ and $0\in\qi(U-U)$.\end{enumerate}}
\end{lemma}

\noindent {\bf Proof.} (ii) Suppose that $x\in\qi(U)$. Then $x\in\qri(U)$ and since
$U-x\subseteq U-U$ and $0\in\qi(U-x)$, the direct implication
follows. The reverse one is a direct consequence of (i) by taking
$V:=\{x\}$. \hfill{$\Box$}

\begin{remark}\label{qi-lp+} {\rm Considering again the setting of
Example \ref{qri-lp+} we get from the second part of the previous
lemma (since $\ell^p_+-\ell^p_+=\ell^p$) that
$$\qi(\ell^p_+)=\qri(\ell^p_+)=\{(x_n)_{n\in\N}\in \ell^p:x_n> 0 \ \forall
n\in\N\}.$$}\end{remark}

\noindent Next we give a useful separation theorem in terms of quasi-relative interior, which will play an important role in the next section when proving the strong duality results.

\begin{theorem}\label{sep-qri-1} (cf. \cite[Theorem 2.7]{qri-siam}) {\it Let $U$ be a nonempty convex subset of $X$ and
$x\in U.$ If $x\not\in\qri(U)$, then there exists $x^*\in X^*,
x^*\neq 0$, such that $$\langle x^*,y\rangle\leq\langle
x^*,x\rangle\mbox{ } \forall y\in U.$$ Viceversa, if there exists
$x^*\in X^*$, $x^*\neq 0$, such that $$\langle x^*,y\rangle\leq\langle x^*,x\rangle\mbox{ } \forall y\in U$$ and $$0\in\qi(U-U),$$ then $x\not\in\qri(U)$.}\end{theorem}

\begin{remark}\label{rem-sep-qri} {\rm (a) The above separation theorem is a
generalization to separated locally convex spaces of a result stated in \cite{Daniele-Giuffre, Daniele-Giuffre-Maugeri} in the framework of real normed spaces (cf. \cite[Remark 2.8]{qri-siam}).

(b) The condition $x\in U$ in Theorem \ref{sep-qri-1} is essential (see
\cite[Remark 2]{Daniele-Giuffre-Maugeri}). However, if $x$ is an arbitrary
element of $X$, an alternative separation theorem has been given by Cammaroto and Di Bella in \cite[Theorem 2.1]{Cammaroto}. Let us mention that some strict separation theorems involving the quasi-relative interior have been provided in \cite{Cammaroto2}.}\end{remark}

\section{General strong duality results}

We briefly recall the general approach for studying conjugate duality by means of the perturbation theory. Consider $X$ and $Y$ separated locally convex spaces and the general optimization problem $$(P_{\Phi})\mbox{ }\mbox{ }\inf_{x\in X}\Phi(x,0),$$  where $\Phi:X\times Y\rightarrow \B$ is the so-called \emph{perturbation function}. By means of the conjugate function of $\Phi$ one can associated to $(P_{\Phi})$ the following dual problem $$(D_{\Phi})\mbox{ }\mbox{ }\sup_{y^*\in Y^*}\{-\Phi^*(0,y^*)\}.$$ The classical Fenchel and Lagrange duality approaches  can be seen as particular instances of this general context. A direct consequence of the Fenchel-Young inequality is the weak duality result, namely $v(D_{\Phi})\leq v(P_{\Phi})$, where $v(P_{\Phi})$ and $v(D_{\Phi})$ are the optimal objective values of the primal and dual problem, respectively. An important issue in conjugate duality is to give conditions which ensure strong duality, the situation when $v(P_{\Phi})=v(D_{\Phi})$ and $(D_{\Phi})$ has an optimal solution. The following regularity conditions were proposed in the literature (see for instance \cite[Theorem 2.7.1]{Zal-carte}):

\begin{center}
\begin{tabular}{r|l}
$(RC^{\Phi}_1)$ \ & \  $\exists x'\in X$ such that $(x',0)\in\dom\Phi$ and $\Phi(x',\cdot)$ is continuous at $0$;
\end{tabular}
\end{center}

\begin{center}
\begin{tabular}{r|l}
$(RC^{\Phi}_2)$ \ & \  $X$ and $Y$ are Fr\'{e}chet spaces, $\Phi$ is lower semicontinuous and\\
\ & \ $0\in\inte(\pr_Y(\dom\Phi))$;
\end{tabular}
\end{center}

\begin{center}
\begin{tabular}{r|l}
$(RC^{\Phi}_3)$ \ & \  $X$ and $Y$ are Fr\'{e}chet spaces, $\Phi$ is lower semicontinuous and\\
\ & \ $0\in\core(\pr_Y(\dom\Phi))$;
\end{tabular}
\end{center}

\begin{center}
\begin{tabular}{r|l}
$(RC^{\Phi}_4)$ \ & \  $X$ and $Y$ are Fr\'{e}chet spaces, $\Phi$ is lower semicontinuous\\
\ & \ $\aff(\pr_Y(\dom\Phi))$ is a closed linear subspace of $Y$ and\\ \ & \ $0\in\icr(\pr_Y(\dom\Phi))$
\end{tabular}
\end{center}

and

\begin{center}
\begin{tabular}{r|l}
$(RC^{\Phi}_5)$ \ & \  $X$ and $Y$ are Fr\'{e}chet spaces, $\Phi$ is lower semicontinuous and\\
\ & \ $0\in\sqri(\pr_Y(\dom\Phi))$.
\end{tabular}
\end{center}

By collecting the corresponding results from \cite{EkTem, Rock-conj-dual, Zal-carte} one can give the following strong duality theorem (see also \cite{b-hab}).

\begin{theorem}\label{strog-dual-class} Let $\Phi:X\times Y\rightarrow\B$ be a proper and convex function. If one of the regularity conditions $(RC^{\Phi}_i)$, $i\in\{1,2,3,4,5\}$, is fulfilled, then $v(P_{\Phi})=v(D_{\Phi})$ and $(D_{\Phi})$ has an optimal solution.\end{theorem}

\begin{remark}\label{impl-class}\rm In case $X$ and $Y$ are Fr\'{e}chet spaces and $\Phi$ is a proper, convex and lower semicontinuous function we have the following relations between the above regularity conditions  $$(RC_1^{\Phi})\Rightarrow(RC_2^{\Phi})\Leftrightarrow(RC_3^{\Phi})\Rightarrow(RC_4^{\Phi})\Leftrightarrow(RC_5^{\Phi}).$$
This fact partially follows from \eqref{int-not-1}, but also by employing the characterization given in \eqref{charicr}. Notice that the \emph{infimal value function} $h : Y \rightarrow \B$, $h(y) = \inf_{x \in X} \Phi(x,y)$, is convex and not necessarily lower semicontinuous, while one has that $\dom h = \pr_Y(\dom\Phi)$. Nevertheless, $\Phi$ being ideally convex, $h$ is a li-convex function (cf. \cite[Proposition 2.2.18]{Zal-carte}). Now by \cite[Theorem 2.2.20]{Zal-carte} it follows that $\core (\dom h) = \inte(\dom h)$, which has as consequence the equivalence of the regularity conditions $(RC_2^{\Phi})$ and $(RC_3^{\Phi})$.
\end{remark}

We introduce in the following some new regularity conditions ensuring strong duality expressed by means of the quasi interior and quasi-relative interior, respectively. We suppose that $0\in\pr_Y(\dom\Phi)$ (which is the same with $v(P_{\Phi})<+\infty$) and $v(P_{\Phi})>-\infty $ (since in case $v(P_{\Phi})=-\infty$ the weak duality result secures strong duality). The function $\Phi$ is suppose to be proper and we assume that the set $\pr_{Y\times\R}(\epi\Phi)$ is convex (this is obviously fulfilled if $\Phi$ is a convex function). Let us notice that this property implies that $\pr_Y(\dom\Phi)$ is convex, too.

Consider the following two regularity conditions:
\begin{center}
\begin{tabular}{r|l}
$(RC^{\Phi}_6)$ \ & \ $0\in\qi\big(\pr_Y(\dom\Phi)\big)$ and\\
& \ $(0,0)\notin\qri\Big[\co\Big(\pr_{Y\times\R}\big(\epi(\Phi-v(P_{\Phi}))\big)\cup\{(0,0)\}\Big)\Big]$
\end{tabular}
\end{center}
and
\begin{center}
\begin{tabular}{r|l}
$(RC^{\Phi}_7)$ \ & \ $0\in\qi\big[\pr_Y(\dom\Phi)-\pr_Y(\dom\Phi)\big]$, $0\in\qri\big(\pr_Y(\dom\Phi)\big)$ and\\
& \ $(0,0)\notin\qri\Big[\co\Big(\pr_{Y\times\R}\big(\epi(\Phi-v(P_{\Phi}))\big)\cup\{(0,0)\}\Big)\Big]$.
\end{tabular}
\end{center}

Notice that
$$\pr\nolimits_{Y\times\R}\big(\epi(\Phi-v(P_{\Phi}))\big)=\pr\nolimits_{Y\times\R}(\epi\Phi)-(0,v(P_{\Phi}))$$
$$=\{(y,r)\in Y\times \R:\exists x\in X\mbox{ such that }\Phi(x,y)-v(P_{\Phi})\leq r\}$$
is also a convex set. Moreover, one can prove that the primal problem has an optimal solution if and only if $(0,0)\in\pr_{Y\times\R}\big(\epi(\Phi-v(P_{\Phi}))\big)$. In the next proposition we study the relations between the two regularity conditions introduced above.

\begin{proposition}\label{comp-qri-phi} Let $\Phi :X \times Y \rightarrow \B$ be a proper function such that $v(P_{\Phi}) \in \R$ and $\pr_{Y\times\R}(\epi\Phi)$ is a convex set. The following statements are true:
\begin{enumerate} \item[(i)] $(RC^{\Phi}_6)\Leftrightarrow(RC^{\Phi}_7)$;
\item[(ii)] if $(P_\Phi)$ has an optimal solution, then $(0,0)\notin\qri\Big[\co\Big(\pr_{Y\times\R}\big(\epi(\Phi-v(P_{\Phi}))\big)\cup\{(0,0)\}\Big)\Big]$ can be equivalently written as $(0,0)\notin\qri\Big(\pr_{Y\times\R}\big(\epi(\Phi-v(P_{\Phi}))\big)\Big)$;
\item[(iii)] if $0\in\qi\big[\pr_Y(\dom\Phi)-\pr_Y(\dom\Phi)\big]$, then $(0,0)\notin\qri\Big[\co\Big(\pr_{Y\times\R}\big(\epi(\Phi-v(P_{\Phi}))\big)\cup\{(0,0)\}\Big)\Big]$ is equivalent to $(0,0)\notin\qi\Big[\co\Big(\pr_{Y\times\R}\big(\epi(\Phi-v(P_{\Phi}))\big)\cup\{(0,0)\}\Big)\Big]$.
\end{enumerate}\end{proposition}

\noindent {\bf Proof.} (i) The equivalence is a direct consequence of Lemma \ref{alte-prop-qri}(ii).

(ii) The statement follows via the comments made before formulating the proposition.

(iii) The direct implication holds trivially, even without the additional assumption. We suppose in the following that $(0,0)\in\qri\Big[\co\Big(\pr_{Y\times\R}\big(\epi(\Phi-v(P_{\Phi}))\big)\cup\{(0,0)\}\Big)\Big]$. Take an arbitrary $(y^*,r^*)\in N_{\co\Big[\pr_{Y\times\R}\big(\epi(\Phi-v(P_{\Phi}))\big)\cup\{(0,0)\}\Big]}(0,0)$. It is enough to show that $(y^*,r^*)=(0,0)$ and the conclusion will follow from Proposition \ref{charact-qi}. By the definition of the normal cone we have
\begin{equation}\label{ineq-phi}\langle y^*,y\rangle+r^*r\leq 0 \ \forall (y,r)\in \pr\nolimits_{Y\times\R}\big(\epi(\Phi-v(P_{\Phi}))\big).\end{equation} Let $x_0\in X$ be such that $\Phi(x_0,0)\in\R$. Taking $(y,r)=(0,\Phi(x_0,0)-v(P_{\Phi})+1)\in \pr\nolimits_{Y\times\R}\big(\epi(\Phi-v(P_{\Phi}))\big)$ in \eqref{ineq-phi} we get $r^*(\Phi(x_0,0)-v(P_{\Phi})+1)\leq 0$, which implies $r^*\leq 0$. Taking into account that $N_{\co\Big[\pr_{Y\times\R}\big(\epi(\Phi-v(P_{\Phi}))\big)\cup\{(0,0)\}\Big]}(0,0)$ is a linear subspace of $Y^*\times\R$ (cf. Proposition \ref{charact-qri}), the same argument applies for $(-y^*,-r^*)$, implying $-r^*\leq 0$. Consequently, $r^*=0$. Using again the inequality in \eqref{ineq-phi} but also the fact that $(-y^*,0)\in N_{\co\Big[\pr_{Y\times\R}\big(\epi(\Phi-v(P_{\Phi}))\big)\cup\{(0,0)\}\Big]}(0,0)$, it follows $$\langle y^*,y\rangle= 0 \ \forall (y,r)\in \pr\nolimits_{Y\times\R}\big(\epi(\Phi-v(P_{\Phi}))\big),$$ which is nothing else than $\langle y^*,y\rangle= 0$ for all $y\in \pr_Y(\dom\Phi)$. Since $y^*$ is linear and continuous, the last relation implies $\langle y^*,y\rangle= 0$ for all $y\in\cl\Big[\cone\big(\pr_Y(\dom\Phi)-\pr_Y(\dom\Phi)\big)\Big]=Y$, hence $y^*=0$ and the conclusion follows.\hfill{$\Box$}

\begin{remark}\label{nice-looking-cond-phi}\rm (a) A sufficient condition which guarantees the relation $0\in\qi\big[\pr_Y(\dom\Phi)-\pr_Y(\dom\Phi)\big]$ in Proposition \ref{comp-qri-phi}(iii) is  $0\in\qi\big(\pr_Y(\dom\Phi)\big)$. This a a consequence of the inclusion  $\pr_Y(\dom\Phi)\subseteq\pr_Y(\dom\Phi)-\pr_Y(\dom\Phi)$.

(b) We have the following implication $$(0,0)\in\qi\Big[\co\Big(\pr\nolimits_{Y\times\R}\big(\epi(\Phi-v(P_{\Phi}))\big)\cup\{(0,0)\}\Big)\Big]\Rightarrow 0\in\qi\big(\pr\nolimits_Y(\dom\Phi)\big).$$ Indeed, suppose that $(0,0)\in\qi\Big[\co\Big(\pr_{Y\times\R}\big(\epi(\Phi-v(P_{\Phi}))\big)\cup\{(0,0)\}\Big)\Big]$. This delivers us the equality   $\cl\Big[\coneco\Big(\pr_{Y\times\R}\big(\epi(\Phi-v(P_{\Phi}))\big)\cup\{(0,0)\}\Big)\Big]=Y\times \R$, hence (cf. \eqref{coneco}) $\cl\Big[\cone\Big(\pr_{Y\times\R}\big(\epi(\Phi-v(P_{\Phi}))\big)\Big)\Big]=Y\times \R$. Since the inclusion $$\cl\Big[\cone\Big(\pr\nolimits_{Y\times\R}\big(\epi(\Phi-v(P_{\Phi}))\big)\Big)\Big]\subseteq \cl\Big(\cone\big(\pr\nolimits_Y(\dom\Phi)\big)\Big)\times\R$$ trivially holds, we get  $\cl\Big(\cone\big(\pr_Y(\dom\Phi)\big)\Big)=Y$, that is $0\in\qi\big(\pr_Y(\dom\Phi)\big)$. Hence the following implication holds $$0\notin\qi\big(\pr\nolimits_Y(\dom\Phi)\big)\Rightarrow (0,0)\notin\qi\Big[\co\Big(\pr\nolimits_{Y\times\R}\big(\epi(\Phi-v(P_{\Phi}))\big)\cup\{(0,0)\}\Big)\Big].$$ Nevertheless, in the regularity conditions $(RC_6^{\Phi})$ and $(RC_7^{\Phi})$ given above one cannot substitute the condition $(0,0)\notin\qi\Big[\co\Big(\pr_{Y\times\R}\big(\epi(\Phi-v(P_{\Phi}))\big)\cup\{(0,0)\}\Big)\Big]$ by the more handleable one  $0\notin\qi\big(\pr_Y(\dom\Phi)\big)$, since the other hypotheses guarantee that $0\in\qi\big(\pr_Y(\dom\Phi)\big)$ (cf. Proposition \ref{comp-qri-phi}(i)).\end{remark}

Let us state now the announced strong duality result.

\begin{theorem}\label{str-dual-phi-qri} Let $\Phi:X\times Y\rightarrow\B$ be a proper function such that $v(P_{\Phi}) \in \R$ and $\pr_{Y\times\R}(\epi\Phi)$ is a convex subset of $Y\times\R$ (the latter is the case if for instance $\Phi$ is a convex function). Suppose that one of the regularity conditions $(RC_i^{\Phi})$, $i\in\{6,7\}$, is fulfilled. Then $v(P_{\Phi})=v(D_{\Phi})$ and $(D_{\Phi})$ has an optimal solution.\end{theorem}

\noindent {\bf Proof.} In view of Proposition \ref{comp-qri-phi}(i) it is enough to give the proof in case $(RC_6^{\Phi})$ is fulfilled, a condition which we assume in the following to be true.

We apply Theorem \ref{sep-qri-1} with $U:=\co\Big(\pr_{Y\times\R}\big(\epi(\Phi-v(P_{\Phi}))\big)\cup\{(0,0)\}\Big)$ and $x:=(0,0)\in U$. Hence there exists $(y^*,r^*)\in Y^*\times\R$, $(y^*,r^*)\neq (0,0)$, such that \begin{equation}\label{sep-phi} \langle y^*,y\rangle + r^* r\leq 0 \ \forall (y,r)\in\pr\nolimits_{Y\times\R}\big(\epi(\Phi-v(P_{\Phi}))\big).\end{equation} We claim that $r^* \leq 0$. Suppose that $r^* >0$. Writing the inequality \eqref{sep-phi} for $(y,r)=(0,\Phi(x_0,0)-v(P_{\Phi})+n)$, where $n\in\N$ is arbitrary and $x_0\in X$ is such that $\Phi(x_0,0)\in\R$, we obtain $r^*(\Phi(x_0,0)-v(P_{\Phi})+n)\leq 0$ for all $n\in\N$. Passing to the limit as $n\rightarrow +\infty$ we get a contradiction. Suppose now that $r^*=0$. Then from \eqref{sep-phi} we obtain $\langle y^*,y\rangle \leq 0$ for all $y\in\pr_Y(\dom\Phi)$, hence $\langle y^*,y\rangle \leq 0$ for all $y\in\cl\Big(\cone\big(\pr_Y(\dom\Phi)\big)\Big)=Y$. The last relation implies $y^*=0$, which contradicts the fact that $(y^*,r^*)\neq (0,0)$. All together, we conclude that $r^* < 0$.

Relation \eqref{sep-phi} ensures $$\left \langle\frac{1}{r^*}y^*,y \right \rangle+r\geq 0 \ \forall (y,r)\in\pr\nolimits_{Y\times\R}\big(\epi(\Phi-v(P_{\Phi}))\big),$$ which is nothing else than $$\left \langle\frac{1}{r^*}y^*,y \right \rangle + \Phi(x,y)-v(P_{\Phi})\geq 0 \ \forall (x,y)\in X\times Y.$$ Further, we have $v(P_{\Phi})\leq\langle (1/r^*) y^*,y\rangle+\Phi(x,y)$, for all $(x,y)\in X\times Y$, which means that $v(P_{\Phi})\leq -\Phi^*(0,-(1/r^*) y^*)\leq v(D_{\Phi})$. As the opposite inequality is always true, we get $v(P_{\Phi})=v(D_{\Phi})$ and $-(1/r^*) y^* \in Y^*$ is an optimal solution of the problem $(D_{\Phi})$.\hfill{$\Box$}

\begin{remark}\label{ren}\rm If we renounce to the condition $(0,0)\notin\qri\Big[\co\Big(\pr_{Y\times\R}\big(\epi(\Phi-v(P_{\Phi}))\big)\cup\{(0,0)\}\Big)\Big]$, then the duality result stated above may fail. See the next sections for several examples in this sense.\end{remark}

In what follows we compare the regularity conditions expressed by means of the quasi interior and quasi-relative interior with the classical ones mentioned at the beginning of this section. We need first an auxiliary result.

\begin{proposition}\label{strong-impl-qi-phi} Suppose that for the primal-dual pair $(P_{\Phi})-(D_{\Phi})$ strong duality holds. Then $(0,0)\notin\qi\Big[\co\Big(\pr_{Y\times\R}\big(\epi(\Phi-v(P_{\Phi}))\big)\cup\{(0,0)\}\Big)\Big]$.\end{proposition}

\noindent {\bf Proof.} By the assumptions we made, there exists $y^*\in Y^*$ such that $v(P_{\Phi})=-\Phi^*(0,y^*)=\inf_{(x,y)\in X\times Y}\{\langle-y^*,y\rangle+\Phi(x,y)\}$. Hence $$\langle-y^*,y\rangle+\Phi(x,y)-v(P_{\Phi})\geq 0 \ \forall (x,y)\in X\times Y,$$ from which we obtain $$\langle-y^*,y\rangle+r\geq 0 \ \forall (y,r)\in \pr\nolimits_{Y\times\R}\big(\epi(\Phi-v(P_{\Phi}))\big),$$ thus $$\langle(y^*,-1),(y,r)\rangle\leq 0 \ \forall (y,r)\in \co\Big(\pr\nolimits_{Y\times\R}\big(\epi(\Phi-v(P_{\Phi}))\big)\cup\{(0,0)\}\Big).$$ The last relation guarantees that $(y^*,-1)\in N_{\co\Big[\pr_{Y\times\R}\big(\epi(\Phi-v(P_{\Phi}))\big)\cup\{(0,0)\}\Big]}(0,0)$, hence we have $N_{\co\Big[\pr_{Y\times\R}\big(\epi(\Phi-v(P_{\Phi}))\big)\cup\{(0,0)\}\Big]}(0,0)\neq\{(0,0)\}$. By applying Proposition \ref{charact-qi} it yields $(0,0)\notin\qi\Big[\co\Big(\pr_{Y\times\R}\big(\epi(\Phi-v(P_{\Phi}))\big)\cup\{(0,0)\}\Big)\Big]$.\hfill{$\Box$}

\begin{proposition}\label{comp-qri-class-phi} Suppose that $X$ and $Y$ are Fr\'{e}chet spaces and $\Phi:X\times Y\rightarrow\B$ is a proper, convex and lower semicontinuous function. The following relations hold $$(RC_1^{\Phi})\Rightarrow(RC_2^{\Phi})\Leftrightarrow(RC_3^{\Phi})\Rightarrow(RC_6^{\Phi})\Leftrightarrow(RC_7^{\Phi}).$$\end{proposition}

\noindent {\bf Proof.} In view of Remark \ref{impl-class} and Proposition \ref{comp-qri-phi}(i) we have to prove only the implication $(RC_3^{\Phi})\Rightarrow(RC_6^{\Phi})$. Let us suppose that $(RC_3^{\Phi})$ is fulfilled. We apply \eqref{int-not-1} and obtain $0\in\qi\big(\pr_Y(\dom\Phi)\big)$. Moreover, the regularity condition $(RC_3^{\Phi})$ ensures strong duality for the pair $(P_{\Phi})-(D_{\Phi})$ (cf. Theorem \ref{strog-dual-class}), hence $(0,0)\notin\qi\Big[\co\Big(\pr_{Y\times\R}\big(\epi(\Phi-v(P_{\Phi}))\big)\cup\{(0,0)\}\Big)\Big]$ (cf.  Proposition \ref{strong-impl-qi-phi}). Applying Proposition \ref{comp-qri-phi}(iii) (see also Remark \ref{nice-looking-cond-phi}(a)) we get that the condition $(RC_6^{\Phi})$ holds and the proof is complete.\hfill{$\Box$}

\begin{remark}\label{rc1-impl-rc6-phi}\rm One can notice that the implications $$(RC_1^{\Phi})\Rightarrow(RC_6^{\Phi})\Leftrightarrow(RC_7^{\Phi})$$ hold in the framework of separated locally convex spaces and for $\Phi:X\times Y\rightarrow\B$ a proper and convex function (nor completeness for the spaces neither lower semicontinuity for the perturbation function is needed here).\end{remark}

\begin{remark}\label{rc4-rc6-nu-phi}\rm In general the conditions $(RC_5^{\Phi})$ and $(RC_6^{\Phi})$ cannot be compared. We underline this fact in the following sections by several examples.\end{remark}

\section{Fenchel duality}

By specializing the results presented in the previous section, we deal in the following with
regularity conditions for the following optimization
problem
$$\mbox{ }(P_F) \ \inf_{x\in X}\{f(x)+g(x)\},$$ where $X$
is a separated locally convex space and $f,g:X\rightarrow\B$ are
proper functions such that $\dom f\cap\dom g\neq\emptyset$.

To this end we consider the perturbation function $\Phi_F:X\times X\rightarrow \B$, defined by $\Phi_F(x,y)=f(x)+g(x-y)$ for all $(x,y)\in X\times X$. The optimal objective value of the primal problem $(P_{\Phi_F})$ is $v(P_{\Phi_F})=v(P_F)$. A simple computation shows that $\Phi_F^*(x^*,y^*)=f^*(x^*+y^*)+g^*(-y^*)$ for all $(x^*,y^*)\in
X^*\times X^*$. The conjugate dual problem to $(P_F)$ obtained by means of the perturbation function $\Phi_F$ is nothing else than its classical \emph{Fenchel dual} problem and it
looks like
$$\mbox{ }(D_F) \ \sup_{y^*\in X^*}\{-f^*(-y^*)-g^*(y^*)\}.$$
One can prove that $\pr_X(\dom\Phi_F)=\dom f-\dom g$, thus the
generalized interior-point regularity conditions
$(RC_i^{\Phi_F})$, $i\in\{1,2,3,4,5\}$, become in this case

\begin{center}
\begin{tabular}{r|l}
$(RC^F_1)$ \ & \ $\exists x'\in\dom f\cap\dom g$ such that $f$ (or $g$) is continuous at $x'$;
\end{tabular}
\end{center}

\begin{center}
\begin{tabular}{r|l}
$(RC^F_2)$ \ & \ $X$ is a Fr\'{e}chet space, $f$ and $g$ are lower semicontinuous and \\ \ & \ $0\in\inte(\dom f-\dom g)$;
\end{tabular}
\end{center}

\begin{center}
\begin{tabular}{r|l}
$(RC^F_3)$ \ & \ $X$ is a Fr\'{e}chet space, $f$ and $g$ are lower semicontinuous and \\ \ & \ $0\in\core(\dom f-\dom g)$;\end{tabular}
\end{center}

\begin{center}
\begin{tabular}{r|l}
$(RC^F_4)$ \ & \ $X$ is a Fr\'{e}chet space, $f$ and $g$ are lower semicontinuous,\\ \ & \ $\aff(\dom f-\dom g)$ is a closed linear subspace of $X$ and \\ \ & \ $0\in\icr(\dom f-\dom g)$\end{tabular}
\end{center}

and

\begin{center}
\begin{tabular}{r|l}
$(RC^F_5)$ \ & \ $X$ is a Fr\'{e}chet space, $f$ and $g$ are lower semicontinuous and \\ \ & \ $0\in\sqri(\dom f-\dom g)$.\end{tabular}
\end{center}

Notice that the condition $(RC_1^{\Phi_F})$ actually asks that $g$
should be continuous at some point $x'\in\dom f\cap\dom g$. But,
when interchanging the roles of $f$ and $g$ one also obtains a
valuable regularity condition for Fenchel strong duality. The
condition $(RC_3^F)$ was considered by Rockafellar (cf.
\cite{Rock-conj-dual}), $(RC_5^F)$ by Rodrigues (cf.
\cite{Rodrigues}), while Gowda and Teboulle proved that $(RC_4^F)$
and $(RC_5^F)$ are equivalent (cf. \cite{Gowda-Teboulle}).

\begin{theorem}\label{strog-dual-class-f} Let $f,g:X\rightarrow\B$ be proper and convex functions. If one of the regularity conditions $(RC^F_i)$, $i\in\{1,2,3,4,5\}$, is fulfilled, then $v(P_F)=v(D_F)$ and $(D_F)$ has an optimal solution.\end{theorem}

\begin{remark}\label{impl-class-f}\rm In case $X$ is a Fr\'{e}chet space and $f,g$ are proper, convex and lower semicontinuous functions we have the following relations between the above regularity conditions (cf. Remark \ref{impl-class}, see also \cite{Gowda-Teboulle, Zal-art} and \cite[Theorem 2.8.7]{Zal-carte}) $$(RC_1^F)\Rightarrow(RC_2^F)\Leftrightarrow(RC_3^F)\Rightarrow(RC_4^F)\Leftrightarrow(RC_5^F).$$ \end{remark}

We suppose in the following that $v(P_F)\in\R$. By some algebraic manipulations we get
$$\pr\nolimits_{X\times\R}\big(\epi(\Phi_F-v(P_F))\big)=\epi f-\epihat(g-v(P_F))$$
$$=\{(x-y,f(x)+g(y)-v(P_F)+\varepsilon):x\in\dom f, y\in\dom g, \varepsilon\geq 0\}.$$
Consider now the following regularity conditions expressed by means of the quasi interior and quasi-relative interior. Besides the conditions $(RC_6^{\Phi_F})$ and $(RC_7^{\Phi_F})$, which in this case are exactly $(RC^F_6)$ and, respectively, $(RC^F_7)$ below, we consider a further one, which we denote by $(RC^F_{6'})$:
\begin{center}
\begin{tabular}{r|l}
$(RC^F_{6'})$ \ & \ $\dom f\cap\qri(\dom
g)\neq\emptyset, \ 0\in\qi(\dom g-\dom g)$ and\\
& \ $(0,0)\notin\qri\Big[\co\Big((\epi
f-\epihat(g-v(P_F)))\cup\{(0,0)\}\Big)\Big]$;
\end{tabular}
\end{center}

\begin{center}
\begin{tabular}{r|l}
$(RC^F_6)$ \ & \ $0\in\qi(\dom f-\dom g)$ and\\
& \ $(0,0)\notin\qri\Big[\co\Big((\epi
f-\epihat(g-v(P_F)))\cup\{(0,0)\}\Big)\Big]$
\end{tabular}
\end{center}

and

\begin{center}
\begin{tabular}{r|l}
$(RC^F_7)$ \ & \ $0\in\qi\big[(\dom f-\dom g)-(\dom f-\dom g)\big]$,  $0\in\qri(\dom f-\dom g)$ and\\
& \ $(0,0)\notin\qri\Big[\co\Big((\epi
f-\epihat(g-v(P_F)))\cup\{(0,0)\}\Big)\Big]$.
\end{tabular}
\end{center}

Let us notice that these three regularity conditions were first introduced in \cite{qri-siam}.

\begin{proposition}\label{comp-qri-f}{\it Let $f,g:X\rightarrow\B$ be proper functions such that $v(P_F) \in \R$ and $\epi f-\epihat(g-v(P_F))$ is a convex subset of $X\times \R$ (the latter is the case if for instance $f$ and $g$ are convex functions). The following statements are true:
\begin{enumerate}\item[(i)] $(RC^F_6)\Leftrightarrow(RC^F_7)$; if, moreover, $f$ and $g$ are convex, then $(RC^F_{6'})\Rightarrow(RC^F_6)\Leftrightarrow(RC^F_7)$;

\item[(ii)] if $(P_F)$ has an optimal solution, then
$(0,0)\notin\qri\Big[\co\Big((\epi f-\epihat(g-v(P_F)))\cup\{(0,0)\}\Big)\Big]$ can be equivalently
written as $(0,0)\notin\qri\Big(\epi f-\epihat(g-v(P_F))\Big)$;

\item[(iii)] if $0\in\qi\big[(\dom f-\dom
g)-(\dom f-\dom g)\big]$, then $(0,0)\notin\qri\Big[\co\Big((\epi
f-\epihat(g-v(P_F)))\cup\{(0,0)\}\Big)\Big]$ is equivalent to
$(0,0)\notin\qi\Big[\co\Big((\epi
f-\epihat(g-v(P_F)))\cup\{(0,0)\}\Big)\Big]$.\end{enumerate}}
\end{proposition}

\noindent {\bf Proof.} In view of Proposition \ref{comp-qri-phi} we have to prove only the second statement of (i). Let us suppose that $f$ and $g$ are convex and $(RC^F_{6'})$ is
fulfilled. By applying Lemma \ref{alte-prop-qri}(i) with $U:=\dom g$ and $V:=\dom f$ we get $0\in\qi(\dom g-\dom f)$, or,
equivalently, $0\in\qi(\dom f-\dom g)$. This means that $(RC^F_6)$ holds.\hfill{$\Box$}

\begin{remark}\label{nice-looking-cond-f}\rm (a) The condition $0\in\qi(\dom f-\dom g)$ implies relation $0\in\qi\big[(\dom f-\dom  g)-(\dom f-\dom g)\big]$ in Proposition \ref{comp-qri-f}(iii).

(b) The following implication holds $$(0,0)\in\qi\Big[\co\Big((\epi
f-\epihat(g-v(P_F)))\cup\{(0,0)\}\Big)\Big]\Rightarrow 0\in\qi(\dom f-\dom g)$$ and in this case similar comments as in Remark \ref{nice-looking-cond-phi}(b) can be made.\end{remark}

A direct consequence of Theorem \ref{str-dual-phi-qri} and
Proposition \ref{comp-qri-f}(i) is the following strong duality
result concerning the pair $(P_F)-(D_F)$. It was first stated in
\cite{qri-siam} under convexity assumptions for the functions
involved.

\begin{theorem}\label{str-dual-qri-f} Let $f,g:X\rightarrow\B$ be proper functions such that $v(P_F) \in \R$ and $\epi f-\epihat(g-v(P_F))$ is a convex subset of $X\times \R$ (the latter is the case if for instance $f$ and $g$ are convex functions). Suppose that either $f$ and $g$ are convex and $(RC^F_{6'})$ is fulfilled, or one of the regularity conditions $(RC^F_i)$, $i\in\{6,7\}$, holds. Then $v(P_F)=v(D_F)$ and $(D_F)$ has an optimal solution.\end{theorem}

When one renounces to the condition $(0,0)\notin\qri\Big[\co\Big((\epi f-\epihat(g-v(P_F)))\cup\{(0,0)\}\Big)\Big]$ the duality
result given above may fail. We give an example (which can be found in \cite{Gowda-Teboulle}, see also \cite{qri-siam}) to show that this assumption is essential.

\begin{example}\label{failure-qri} {\rm Consider
the real Hilbert space $X=\ell^2(\N)$ and the sets
$$C=\{(x_n)_{n\in\N}\in \ell^2:x_{2n-1}+x_{2n}=0 \ \forall n\in\N\}$$ and
$$S=\{(x_n)_{n\in\N}\in \ell^2:x_{2n}+x_{2n+1}=0 \ \forall n\in\N\},$$ which
are closed linear subspaces of $\ell^2$ and satisfy $C\cap
S=\{0\}$. Define the functions $f,g:\ell^2\rightarrow\B$ by
$f=\delta_C$ and, respectively, $g(x)=x_1+\delta_S(x)$ for all
$x=(x_n)_{n\in\N}\in\ell^2$. One can see that $f$ and $g$ are
proper, convex and lower semicontinuous functions with $\dom f=C$
and $\dom g=S$. As $v(P_F)=0$ and $v(D_F)=-\infty$ (cf.
\cite[Example 3.3]{Gowda-Teboulle}), we have a duality gap between
the optimal objective values of the primal problem and its Fenchel
dual problem. Moreover, $S-C$ is dense in $\ell^2$ (cf.
\cite{Gowda-Teboulle}), thus $\cl\big(\cone(\dom f-\dom
g)\big)=\cl(C-S)=\ell^2$. The last relation implies $0\in\qi(\dom
f-\dom g)$. Notice that in \cite[Remark 3.12(b)]{qri-siam} it has
been proved that $(0,0)\in\qri\Big(\epi
f-\epihat(g-v(P_F))\Big)$.}\end{example}

Let us give in the following an example which illustrates the
applicability of the strong duality result introduced above. It
has been considered in \cite[Example 3.13]{qri-siam}, however we
give here the details for the sake of completeness.

\begin{example}\label{ex-qri-f} {\rm Consider the real Hilbert space $\ell^2=\ell^2(\N)$. We define the functions
$f,g:\ell^2\rightarrow\B$ by
$$f(x)=\left\{
\begin{array}{ll}
\|x\|, & \mbox {if } x\in x^0-\ell^2_+,\\
+\infty, & \mbox{otherwise}
\end{array}\right.$$ and $$g(x)=\left\{
\begin{array}{ll}
\langle c,x\rangle, & \mbox {if } x\in \ell^2_+,\\
+\infty, & \mbox{otherwise},
\end{array}\right.$$ respectively, where $x^0, c\in\ell^2_+$ are arbitrarily chosen such that $x^0_n>0$ for all $n\in\N$. Note that $$v(P_F)=\inf\limits_{x\in
\ell^2_+\cap(x^0-\ell^2_+)}\{\|x\|+\langle c,x\rangle\}=0$$ and the infimum is attained at $x=0$. We
have $\dom f=x^0-\ell^2_+=\{(x_n)_{n\in\N}\in \ell^2:x_n\leq
x_n^0 \ \forall n\in\N\}$ and $\dom g=\ell^2_+$. By using
Example \ref{qri-lp+} we get $$\dom f\cap\qri(\dom
g)=\{(x_n)_{n\in\N}\in \ell^2:0<x_n\leq x_n^0 \ \forall
n\in\N\}\neq\emptyset.$$ Also, $\cl\big(\cone(\dom g-\dom
g)\big)=\ell^2$ and so $0\in\qi(\dom g-\dom g)$. Further,  $$\epi
f-\epihat(g-v(P_F))=\{(x-y,\|x\|+\langle
c,y\rangle+\varepsilon):x\in x^0-\ell^2_+,y\in
\ell^2_+,\varepsilon\geq 0\}.$$

In the following we prove that $(0,0)\notin\qri\Big(\epi
f-\epihat(g-v(P_F))\Big)$. Assuming the contrary, one would have
that the set $\cl\Big[\cone\Big(\epi
f-\epihat(g-v(P_F))\Big)\Big]$ is a linear subspace of
$\ell^2\times\R$. Since $(0, 1) \in \cl\Big[\cone\Big(\epi
f-\epihat(g-v(P_F))\Big)\Big]$ (take $x=y=0$ and $\varepsilon
=1$),  $(0,-1)$ must belong to this set, too. On the
other hand, one can easily see that for all $(x,r)$ belonging to
$\cl\Big[\cone\Big(\epi f-\epihat(g-v(P_F))\Big)\Big]$ it holds $r
\geq 0$. This leads to the desired contradiction.

Hence the regularity condition $(RC^F_{6'})$ is fulfilled, thus
strong duality holds (cf. Theorem \ref{str-dual-qri-f}). On the
other hand, $\ell^2$ is a Fr\'{e}chet space (being a Hilbert
space), the functions $f$ and $g$ are proper, convex and lower
semicontinuous and, as $\sqri(\dom f-\dom
g)=\sqri(x^0-\ell^2_+)=\emptyset$, none of the regularity
conditions $(RC_i^F)$, $i\in\{1,2,3,4,5\}$, listed at the
beginning of this section can be applied for this optimization
problem.

As for all $x^* \in \ell^2$ it holds
$g^*(x^*)=\delta_{c-\ell^2_+}(x^*)$ and (cf. \cite[Theorem
2.8.7]{Zal-carte})
$$f^*(-x^*) = \inf_{x_1^* + x_2^* = -x^*} \{\| \cdot \|^*(x_1^*) + \delta^*_{x^0 - l^2_+}(x_2^*) \}
= \inf_{\substack{x_1^* + x_2^* = -x^*,\\ \|x_1^*\| \leq 1, x_2^*
\in \ell^2_+}} \langle x_2^*, x^0 \rangle ,$$  the optimal
objective value of the Fenchel dual problem is
$$v(D_F) = \sup_{\substack{x_2^* \in \ell^2_+ - c - x_1^*,\\ \|x_1^*\| \leq 1,
x_2^* \in \ell^2_+}}\langle - x_2^*, x^0 \rangle  = \sup_{x_2^*
\in \ell^2_+}\langle - x_2^*, x^0 \rangle  = 0$$ and $x_2^*=0$ is
the optimal solution of the dual.}\end{example}

The following example (see also \cite[Example 2.5]{teza-csetnek})
underlines the fact that in general the regularity condition
$(RC_6^F)$ (and automatically also $(RC_7^F)$) is weaker than
$(RC_{6'}^F)$ (see also Example \ref{nu-cl-qri-da} below).

\begin{example}\label{ex-rc7-rc6-f} \rm Consider the real Hilbert space $\ell^2(\R)$ and the functions $f,g:\ell^2(\R)\rightarrow \B$ defined for all $s\in\ell^2(\R)$ by $$f(s)=\left\{
\begin{array}{ll}
s(1), & \mbox {if } s\in \ell^2_+(\R),\\
+\infty, & \mbox{otherwise}
\end{array}\right.$$ and $$g(s)=\left\{
\begin{array}{ll}
s(2), & \mbox {if } s\in \ell^2_+(\R),\\
+\infty, & \mbox{otherwise},
\end{array}\right.$$ respectively. The optimal objective value of the primal problem is
$$v(P_F)=\inf_{s\in\ell^2_+(\R)}\{s(1)+s(2)\}=0$$ and $s=0$ is an optimal solution
(let us notice that $(P_F)$ has infinitely many optimal
solutions). We have $\qri(\dom
g)=\qri(\ell^2_+(\R))=\emptyset$ (cf. Example \ref{qri-vid}), hence
the condition $(RC_{6'}^F)$ fails. In the following we show that
$(RC_6^F)$ is fulfilled. One can prove that $\dom f-\dom
g=\ell^2_+(\R)-\ell^2_+(\R)=\ell^2(\R)$, thus $0\in\qi(\dom f-\dom
g)$. Like in the previous example, we have
$$\epi f-\epihat(g-v(P_F))=\{(s-s',s(1)+s'(2)+\varepsilon):s,s'\in\ell^2_+(\R),\varepsilon\geq
0\}$$ and with the same technique one can show that $(0,0)\notin\qri\Big(\epi
f-\epihat(g-v(P_F))\Big)$, hence the condition $(RC_6^F)$ holds.

Let us take a look at the formulation of the dual problem. To this
end we have to calculate the conjugates of $f$ and $g$. Let us
recall that the scalar product on $\ell^2(\R)$,
$\langle\cdot,\cdot\rangle:\ell^2(\R)\times \ell^2(\R)\rightarrow
\R$ is defined by $\langle s,s'\rangle=\sup_{F\subseteq \R,
F\mbox{\tiny finite}}\sum_{r\in F}s(r)s'(r)$, for
$s,s'\in\ell^2(\R)$ and that the dual space
$\big(\ell^2(\R)\big)^*$ is identified with $\ell^2(\R)$. For an
arbitrary $u\in \ell^2(\R)$ we have
$$f^*(u)=\sup_{s\in\ell^2_+(\R)}\{\langle u,s\rangle-s(1)\}
=\sup_{s\in\ell^2_+(\R)}\left\{\sup\limits_{F\subseteq \R,
F\mbox{\tiny finite}}\sum_{r\in F}u(r)s(r)-s(1)\right\}$$$$
=\sup\limits_{F\subseteq \R, F\mbox{\tiny
finite}}\left\{\sup_{s\in\ell^2_+(\R)}\Big\{\sum_{r\in
F}u(r)s(r)-s(1)\Big\}\right\}.$$ If there exists $r \in \R
\setminus \{1\}$ with $u(r) > 0$ or if $u(1) > 1$, then one has
$f^*(u) = +\infty$. Assuming the contrary, for every finite subset
$F$ of $\R$, independently from the fact that $1$ belongs to $F$
or not, it holds $\sup_{s\in\ell^2_+(\R)} \{\sum_{r\in
F}u(r)s(r)-s(1)\} = 0$. Consequently,
$$f^*(u)=\left\{
\begin{array}{ll}
0, & \mbox {if } u(r)\leq 0 \ \forall r\in\R\setminus\{1\}\mbox{ and }u(1)\leq 1,\\
+\infty, & \mbox{otherwise}.
\end{array}\right.$$ Similarly one can provide a formula for $g^*$ and  in this way we obtain that $v(D_F)=0$
and that $u=0$ is an optimal solution of the dual (in fact $(D_F)$ has infinitely many optimal solutions).
\end{example}

From Proposition \ref{strong-impl-qi-phi} we obtain the next result.

\begin{proposition}\label{strong-impl-qi-f} Suppose that for the primal-dual pair $(P_F)-(D_F)$ strong duality holds. Then $(0,0)\notin\qi\Big[\co\Big((\epi f-\epihat(g-v(P_F)))\cup\{(0,0)\}\Big)\Big]$. \end{proposition}

A comparison of the above regularity conditions is provided in the following.

\begin{proposition}\label{comp-qri-class-f} Suppose that $X$ is a Fr\'{e}chet space and $f,g:X\rightarrow\B$ are proper, convex and lower semicontinuous functions. The following relations hold $$(RC_1^F)\Rightarrow(RC_2^F)\Leftrightarrow(RC_3^F)\Rightarrow(RC_6^F)\Leftrightarrow(RC_7^F).$$\end{proposition}

\begin{remark}\label{rc1-impl-rc6-f}\rm One can notice that the implications $$(RC_1^F)\Rightarrow(RC_6^F)\Leftrightarrow(RC_7^F)$$ hold in the framework of separated locally convex spaces and for $f,g : X \rightarrow\B$ proper and convex functions (nor completeness for the space neither lower semicontinuity for the functions is needed here).\end{remark}

In general the conditions $(RC_i^F)$, $i\in\{4,5\}$, cannot be compared with $(RC_i^F)$, $i\in\{6',6,7\}$. Example \ref{ex-qri-f} provides a situation when $(RC_i^F)$, $i\in\{6',6,7\}$, are fulfilled, unlike $(RC_i^F)$, $i\in\{4,5\}$. In the following example the conditions $(RC_i^F)$, $i\in\{4,5\}$, are fulfilled, while $(RC_i^F)$, $i\in\{6',6,7\}$, fail.

\begin{example}\label{sqri-qri-nu}\rm Consider $(X,\|\cdot\|)$ a nonzero real Banach space, $x_0^*\in X^*\setminus\{0\}$ and the functions $f,g:X\rightarrow\B$ defined by $f=\delta_{\Ker x_0^*}$ and $g=\|\cdot\|+\delta_{\Ker x_0^*}$, respectively. The optimal objective value of the primal problem is $$v(P_F)=\inf_{x\in\Ker x_0^*}\|x\|=0$$ and $\bar x=0$ is the unique optimal solution of $(P_F)$. The functions $f$ and $g$ are proper, convex and lower semicontinuous. Further, $\dom f-\dom g=\Ker x_0^*$, which is a closed linear subspace of $X$, hence $(RC_i^F)$, $i\in\{4,5\}$, are fulfilled. Moreover, $\dom g-\dom g=\dom f-\dom g=\Ker x_0^*$ and it holds $\cl(\Ker x_0^*)=\Ker x_0^*\neq X$.
Thus $0\notin\qi(\dom g-\dom g)$ and $0\notin\qi(\dom f-\dom g)$
and this means that all the three regularity conditions
$(RC_i^F)$, $i\in\{6',6,7\}$, fail.

The conjugate functions of $f$ and $g$ are  $f^*=\delta_{(\Ker
x_0^*)^\bot}=\delta_{\R x_0^*}$ and
$g^*=\delta_{B_*(0,1)}\Box\delta_{\R x_0^*}=\delta_{B_*(0,1)+\R
x_0^*}$ (cf. \cite[Theorem 2.8.7]{Zal-carte}), respectively, where
$B_*(0,1)$ is the closed unit ball of the dual space $X^*$. Hence
$v(D_F)=0$ and the set of optimal solutions of $(D_F)$ coincides
with $\R x_0^*$. Finally, let us notice that instead of $\Ker
x_0^*$ one can consider any closed linear subspace $S$ of $X$ such
that $S\neq X$.\end{example}

Let us mention that besides the above mentioned generalized
interior-point regularity conditions, one can meet in the
literature so-called \emph{closedness-type regularity conditions}.
For the primal-dual pair $(P_F)-(D_F)$ this condition has been
first considered by Burachik and Jeyakumar in Banach spaces (cf.
\cite{bur-jeyak}) and by Bo\c t and Wanka in separated locally
convex spaces (cf. \cite{BW}) and it looks like:

\begin{center}
\begin{tabular}{r|l}
$(RC^F_8)$ \ & \ $f$ and $g$ are lower semicontinuous and\\
           & \ $\epi f^*+\epi g^*$ is closed in $(X^*,w(X^*,X))\times\R$.
\end{tabular}
\end{center}

We have the following duality result (cf. \cite{BW}).

\begin{theorem}\label{closed-cond-f} Let $f,g:X\rightarrow \B$ be proper and convex functions such that $\dom f\cap\dom g\neq\emptyset$. If $(RC^F_8)$ is fulfilled, then \begin{equation}\label{stab-str}(f+g)^*(x^*)=\min\{f^*(x^*-y^*)+g^*(y^*):y^*\in X^*\} \ \forall x^*\in X^*.\end{equation}\end{theorem}

\begin{remark}\label{obs-stab-str}\rm (a) Let us notice that condition \eqref{stab-str} is referred in the literature as \emph{stable strong duality} (see \cite{b-hab,BuJeyakW,simons} for more details) and obviously guarantees strong duality for $(P_F)-(D_F)$. When $f,g:X\rightarrow \B$ are proper, convex and lower semicontinuous functions with $\dom f\cap\dom g\neq\emptyset$ one has in fact that  $(RC^F_8)$ is fulfilled if and only if \eqref{stab-str}
holds (cf. \cite[Theorem 3.2]{BW}).

(b) In case $X$ is a Fr\'{e}chet space and $f,g$ are proper, convex and lower semicontinuous functions we have the following relations between the regularity conditions considered for the primal-dual pair $(P_F)-(D_F)$ (cf. \cite{BW}, see also \cite{Gowda-Teboulle} and \cite[Theorem 2.8.7]{Zal-carte}) $$(RC_1^F)\Rightarrow(RC_2^F)\Leftrightarrow(RC_3^F)\Rightarrow(RC_4^F)\Leftrightarrow(RC_5^F)\Rightarrow(RC_8^F).$$ We refer to \cite{BW,b-hab,simons,bur-jeyak} for several examples showing that in general the implications above are strict. The implication $(RC_1^F)\Rightarrow(RC_8^F)$ holds in the general setting of separated locally convex spaces (in the hypotheses that $f,g$ are proper, convex and lower semicontinuous).\end{remark}

We observe that if $X$ is a finite-dimensional space and $f,g$ are proper, convex and lower semicontinuous, then $(RC_{6'}^F)\Rightarrow(RC_6^F)\Leftrightarrow(RC_7^F)\Rightarrow(RC_8^F)$. However, in the infinite-dimensional setting this is no longer true. In the following two examples the condition $(RC_8^F)$ is fulfilled, unlike $(RC_i^F)$, $i\in\{6',6,7\}$ (we refer to \cite{BW,b-hab,simons,bur-jeyak,li-fang-lopez} for examples in the finite-dimensional setting).

\begin{example}\label{cl-nu-qri-f-1}\rm Consider the same setting as in Example \ref{sqri-qri-nu}. We know that $(RC_5^F)$ is fulfilled, hence also $(RC_8^F)$ (cf. Remark \ref{obs-stab-str}(b)). This is not surprising, since $\epi f^*+\epi g^*=(B_*(0,1)+\R x_0^*)\times [0,\infty)$, which is closed in $(X^*,w(X^*,X))\times\R$ (by the Banach-Alaoglu Theorem the unit ball $B_*(0,1)$ is compact in $(X^*,w(X^*,X))$). As shown in Example \ref{sqri-qri-nu}, none of the regularity conditions $(RC_i^F)$, $i\in\{6',6,7\}$ is fulfilled.\end{example}

\begin{example}\label{cl-nu-qri-f-2}\rm Consider the real Hilbert space $\ell^2(\R)$ and the functions $f,g:\ell^2(\R)\rightarrow\B$, defined by $f=\delta_{\ell^2_+(\R)}$ and $g=\delta_{-\ell^2_+(\R)}$, respectively. We have $\qri(\dom f-\dom g)=\qri\big(\ell^2_+(\R)\big)=\emptyset$ (cf. Example \ref{qri-vid}), hence all the generalized interior-point regularity conditions $(RC_i^F)$, $i\in\{1,2,3,4,5,6',6,7\}$ fail (see also Proposition \ref{comp-qri-f}(i)). The conjugate functions of $f$ and $g$ are $f^*=\delta_{-\ell^2_+(\R)}$ and $g^*=\delta_{\ell^2_+(\R)}$, respectively, hence $\epi f^*+\epi g^*=\ell^2(\R)\times[0,\infty)$, that is the condition $(RC_8^F)$ holds. One can see that $v(P_F)=v(D_F)=0$ and  $y^*=0$ is an optimal solution of the dual problem.\end{example}

The next question that we address concerns the relation between
the generalized interior-point condition $(RC_6^F)$ and the
closedness-type one $(RC_8^F)$. In Example \ref{cl-nu-qri-f-2}
(see also Example \ref{cl-nu-qri-f-1}) we have a situation when
the second is fulfilled, while the first one fails. In the
following we provide two examples for which this time $(RC_6^F)$
is fulfilled, unlike $(RC_8^F)$. In this way we give a negative
answer to an open problem stated in \cite{li-fang-lopez},
concomitantly proving that in general $(RC_6^F)$ (and
automatically also $(RC_7^F)$) and $(RC_8^F)$ are not comparable.

\begin{example}\label{nu-cl-qri-da}\rm  Like in Example \ref{failure-qri},
consider the real Hilbert space $X=\ell^2(\N)$ and the sets
$$C=\{(x_n)_{n\in\N}\in \ell^2:x_{2n-1}+x_{2n}=0 \ \forall n\in\N\}$$ and
$$S=\{(x_n)_{n\in\N}\in \ell^2:x_{2n}+x_{2n+1}=0 \ \forall n\in\N\},$$ which
are closed linear subspaces of $\ell^2$ and satisfy $C\cap
S=\{0\}$. Define the functions $f,g:\ell^2\rightarrow\B$ by
$f=\delta_C$ and $g=\delta_S$, respectively, which are proper,
convex and lower semicontinuous. The optimal objective value of
the primal problem is $v(P_F)=0$ and $\overline{x}=0$ is the
unique optimal solution of $v(P_F)$. Moreover, $S-C$ is dense in
$\ell^2$ (cf. \cite[Example 3.3]{Gowda-Teboulle}), thus
$\cl\big(\cone(\dom f-\dom g)\big)=\cl(C-S)=\ell^2$. This implies
$0\in\qi(\dom f-\dom g)$. Further, one has
$$\epi f-\epihat(g-v(P_F))=\{(x-y,\varepsilon):x\in C,y\in
S,\varepsilon\geq 0\}=(C-S)\times[0,+\infty)$$ and
$\cl\Big[\cone\Big(\epi f-\epihat(g-v(P_F))\Big)\Big]=\ell^2\times
[0,+\infty)$, which is not a linear subspace of $\ell^2\times\R$,
hence $(0,0)\notin\qri\Big(\epi f-\epihat(g-v(P_F))\Big)$. All
together, we get that the condition $(RC_6^F)$ is fulfilled, hence
strong duality holds (cf. Theorem \ref{str-dual-qri-f}). One can
prove that $f^*=\delta_{C^\bot}$ and $g^*=\delta_{S^\bot}$, where
$$C^\bot=\{(x_n)_{n\in\N}\in \ell^2:x_{2n-1}=x_{2n} \ \forall
n\in\N\}$$ and
$$S^\bot=\{(x_n)_{n\in\N}\in \ell^2:x_1=0, x_{2n}=x_{2n+1} \ \forall n\in\N\}.$$
Further, $v(D_F)=0$ and the set of optimal solutions of the dual
problem is exactly $C^\bot\cap S^\bot=\{0\}$.

We show that $(RC_8^F)$ is not fulfilled. Let us consider the
element $e^1\in\ell^2$, defined by $e^1_1=1$ and $e^1_k=0$ for all
$k\in\N\setminus\{1\}.$ We compute
$(f+g)^*(e^1)=\sup_{x\in\ell^2}\{\langle
e^1,x\rangle-f(x)-g(x)\}=0$ and $(f^*\Box
g^*)(e^1)=\delta_{C^\bot+S^\bot}(e^1)$. If we suppose that $e^1\in
C^\bot+S^\bot$, then we would have $(e^1+S^\bot)\cap
C^\bot\neq\emptyset$. However, it has been proved in \cite[Example
3.3]{Gowda-Teboulle} that $(e^1+S^\bot)\cap C^\bot=\emptyset$.
This shows that $(f^*\Box g^*)(e^1)=+\infty>0=(f+g)^*(e^1)$. Via
Theorem \ref{closed-cond-f} follows that the condition $(RC_8^F)$
is not fulfilled and, consequently, $(RC_i^F)$,
$i\in\{1,2,3,4,5\}$, fail, too (cf. Remark \ref{obs-stab-str}(b)),
unlike condition $(RC_6^F)$. Looking at $(RC_{6'}^F)$, one can see
that this condition is also not fulfilled, since  $0\in\qi(\dom
g-\dom g)$ does not hold.

Finally, let us notice that one can prove directly that $(RC_8^F)$
is not fulfilled. Indeed, we have $\epi f^*+\epi
g^*=(C^\bot+S^\bot)\times [0,\infty)$. As in \cite[Example
3.3]{Gowda-Teboulle}, one can show that $C^\bot+S^\bot$ is dense
in $\ell^2$. If we suppose that $C^\bot+S^\bot$ is closed, we
would have $C^\bot+S^\bot=\ell^2$, which is a contradiction, since
$e^1\notin C^\bot+S^\bot$.
\end{example}

\begin{example}\label{nu-cl-qri-da2}\rm The example we consider in the following is inspired by \cite[Example 11.3]{simons}. Consider $X$ an arbitrary Banach space, $C$ a convex and closed subset of $X$ and $x_0$ an extreme point of $C$ which is not a support point of $C$. Taking for instance $X=\ell^2$, $1<p<2$ and $C:=\left \{x \in \ell^2: \sum_{n=1}^{\infty}|x_n|^p \leq 1 \right \}$ one can find extreme points in $C$ that are not support points (see \cite{simons}). Consider the functions $f,g :X \rightarrow \overline \R$ defined as $f=\delta_{x_0-C}$ and $g=\delta_{C-x_0}$, respectively. They are both proper, convex and lower semicontinuous and fulfill, as $x_0$ is an extreme point of $C$, $f+g= \delta_{\{0\}}$. Thus $v(P_F) = 0$ and $\overline{x}=0$ is the unique optimal solution of $(P_F)$. We show that, different to the previous example, $(RC^F_{6'})$ is fulfilled and this will guarantee that both $(RC^F_{6})$ and $(RC^F_{7})$ are valid, too. To this end we notice first that $x_0 \in \qi (C)$.
Assuming the contrary, one would have that there exists $x^* \in
X^* \setminus \{0\}$ such that $\langle x^*,x_0 \rangle = \sup_{x
\in C} \langle x^*,x \rangle$ (cf. Proposition \ref{charact-qi}),
contradicting the hypothesis that $x_0$ is not a support point of
$C$. This means that $x_0 \in \qri(C)$, too, and so $0 \in \dom f
\cap \qri(\dom g)$. Further, since it holds $\cl(\cone(C-x_0))
\subseteq \cl(\cone(C-C))$, we have $\cl(\cone(C-C)) = X$ and from
here $0 \in \qi(C-C) = \qi(\dom g - \dom g)$. Noticing that
$$\epi f-\epihat(g-v(P_F))=\{(x-y,\varepsilon):x,y \in C,\varepsilon\geq 0\}=(C-C)\times[0,+\infty),$$
it follows that $\cl\Big[\cone\Big(\epi
f-\epihat(g-v(P_F))\Big)\Big]=X\times[0,+\infty)$, which is not a
linear subspace of $X\times\R$. Thus $(0,0)\notin\qri\Big(\epi
f-\epihat(g-v(P_F))\Big)$ and this has as consequence the fact
that $(RC^F_{6'})$ is fulfilled. Hence strong duality holds (cf.
Theorem \ref{str-dual-qri-f}), $v(D_F)=0$ and $0$ is an optimal
solution of the dual problem.

We show that $(RC_8^F)$ is not fulfilled. Assuming the contrary, one would have that the equality in \eqref{stab-str} holds for all $x^* \in X^*$. On the other hand, in \cite[Example 11.3]{simons} it is proven that this is the case only when $x^*=0$ and this provides the desired contradiction.
\end{example}

\begin{remark}\label{qri-f-g-a}\rm Consider the following optimization
problem
$$\mbox{ }(P_F^A) \ \inf_{x\in X}\{f(x)+(g\circ A)(x)\},$$ where $X$
and $Y$ are separated locally convex spaces with topological
dual spaces $X^*$ and $Y^*$, respectively, $A:X\rightarrow Y$ is a
linear continuous mapping, $f:X\rightarrow\B$ and
$g:Y\rightarrow\B$ are proper functions such that $A(\dom
f)\cap\dom g\neq\emptyset$. The Fenchel dual problem to $(P_F^A)$ is
$$\mbox{ }(D_F^A) \ \sup_{y^*\in Y^*}\{-f^*(-A^*y^*)-g^*(y^*)\}.$$ We denote by
$v(P_F^A)$ and $v(D_F^A)$ the optimal objective values of the
primal and the dual problem, respectively, and suppose that
$v(P_F^A)\in\R$. We consider the set $$A\times{\id}_{\R}(\epi
f)=\{(Ax,r)\in Y\times\R:f(x)\leq r\}.$$ By using the approach
presented in the previous section one can provide similar
discussions regarding strong duality for the primal-dual pair
$(P_F^A)-(D_F^A)$. To this end, one has to define the perturbation
function $\Phi_A:X\times Y\rightarrow \B$,
$\Phi_A(x,y)=f(x)+g(Ax-y)$ for all $(x,y)\in X\times Y$. Notice
that in this case $\pr_Y(\dom\Phi_A)=A(\dom f)-\dom g$ and
$\pr_{Y\times\R}\big(\epi(\Phi_A-v(P_F^A))\big)=A\times{\id}_{\R}(\epi
f)-\epihat(g-v(P_F^A))$. We refer also to \cite{qri-siam} where
the strong duality results for the pair $(P_F^A)-(D_F^A)$ are
deduced from the corresponding ones given for $(P_F)-(D_F)$.

Let us notice that Borwein and Lewis gave in
\cite{Borwein-Lewis} some regularity conditions by means of the
quasi-relative interior, in order to guarantee strong duality for $(P_F^A)$ and $(D_F^A)$. However, they considered a more restrictive case, namely when the codomain of the operator $A$
is finite-dimensional. Here we have considered the more general case, when both  spaces $X$ and $Y$ are infinite-dimensional.\end{remark}

\section{Lagrange duality}

The aim of this section is to particularize the investigations made in section 4 to another classical duality concept, the \emph{Lagrange duality}, this time in connection
to the optimization problem with geometric and cone constraints
$$\begin{array}{rl}
\mbox{ }(P_L) \  & \inf\limits_{x \in T} f(x),\\
& T = \{x \in S : g(x) \in -C\}
\end{array}$$
where $X$ and $Z$ are two separated locally convex spaces, the
latter being partially ordered by the nonempty convex cone $C
\subseteq Z$, $S \subseteq X$ is a nonempty set and $f : X
\rightarrow \overline \R$ and $g : X \rightarrow Z^{\bullet}$ are
proper functions fulfilling $\dom f \cap S \cap g^{-1}(-C) \neq
\emptyset$.

Considering as perturbation function $\Phi_L :X \times Z
\rightarrow \B$ defined by $\Phi_L(x,z)=f(x)+ \delta_{\{y \in S :
g(y) \in z-C\}}(x)$ for all $(x,z)\in X\times Y$, the optimal
objective value of the primal problem $(P_{\Phi_L})$ is nothing
else than $v(P_{\Phi_L})=v(P_L)$. The conjugate of $\Phi_L$ is
$\Phi_L^*(x^*,z^*)= (f+ (-z^*g) + \delta_S)^*(x^*) +
\delta_{-C^*}(z^*)$ for all $(x^*,z^*)\in X^*\times Y^*$ and so
the conjugate dual problem to $(P_L)$ obtained by means of
$\Phi_L$ is nothing else than its classical \emph{Lagrange dual}
problem. This looks like
$$\mbox{ }(D_L) \ \sup_{z^*\in C^*} \inf_{x \in S}\{f(x) + \langle z^*,g(x)\rangle\}.$$

The regularity condition $(RC_1^{\Phi_L})$ states in this
particular case that there exists $x' \in \dom f \cap S \cap
g^{-1}(-C)$ such that the function $z \mapsto f(x') + \delta_S(x')
+ \delta_{g(x')+C}(z)$ is continuous at $0$, which is the same
with saying that there exists $x' \in \dom f \cap S \cap
g^{-1}(-C)$ such that $0 \in \inte(g(x') + C)$ or, equivalently,
with asking that
\begin{center}
\begin{tabular}{r|l}
$(RC^{L}_{1})$ \ & \ $\exists x' \in \dom f \cap S$ such that
$g(x') \in -\inte(C)$.
\end{tabular}
\end{center}
This is nothing else than the classical \emph{Slater constraint
qualification}.

We come now to the class of regularity conditions which assume
that $X$ and $Z$ are Fr\' echet spaces. One has $\pr_Z(\dom
\Phi_L) = g(\dom f \cap S \cap \dom g) + C$ and in order to
guarantee the lower semicontinuity for $\Phi_L$ it is enough to
assume that $S$ is closed, $f$ is lower semicontinuous and $g$ is $C$-epi closed. Under these assumptions the epigraph of the
perturbation function
$$\epi \Phi_L = \{(x,z,r) \in  X \times Z \times \R : (x,r) \in \epi f\} \cap (S \times Z \times \R) \cap (\epi\nolimits_C g \times \R)$$
is a closed set. These lead to the following
regularity conditions (cf. \cite{b-hab})

\begin{center}
\begin{tabular}{r|l}
$(RC^{L}_{2})$ & \ $X$ and $Z$ are Fr\' echet spaces, \!$S$ is
closed, \!$f$ is lower semicontinuous,\\
& \ $g$ is $C$-epi closed and $0 \in \inte \big (g(\dom f \cap S
\cap \dom g) + C \big)$;
\end{tabular}
\end{center}

\begin{center}
\begin{tabular}{r|l}
$(RC^{L}_{3})$ & \ $X$ and $Z$ are Fr\' echet spaces, \!$S$ is
closed, \!$f$ is lower semicontinuous,\\
& \ $g$ is $C$-epi closed and $0 \in \core \big (g(\dom f \cap S
\cap \dom g) + C \big)$;
\end{tabular}
\end{center}

\begin{center}
\begin{tabular}{r|l}
$(RC^{L}_4)$ & \ $X$ and $Z$ are Fr\' echet spaces, \!$S$ is
closed, \!$f$ is lower semicontinuous,\\
& \ $g$ is $C$-epi closed, $\aff \big (g(\dom f \cap S
\cap \dom g) + C \big)$ is a closed linear subspace\\
& \ and $0 \in \icr \big (g(\dom f \cap S \cap \dom g) + C \big)$
\end{tabular}
\end{center}

and

\begin{center}
\begin{tabular}{r|l}
$(RC^{L}_{5})$ & \ $X$ and $Z$ are Fr\' echet spaces, \!$S$ is
closed, \!$f$ is lower semicontinuous,\\
& \ $g$ is $C$-epi closed, $0 \in \sqri \big (g(\dom f \cap S \cap
\dom g) + C \big)$.
\end{tabular}
\end{center}

The condition $(RC_3^L)$ was considered in Banach spaces by
Rockafellar in \cite{Rock-conj-dual}, while a particular
formulation of $(RC^L_{5})$ has been stated for linear programming
problems in \cite{Jeyak-Wolk}. Let us notice that some stronger
versions of the above regularity conditions have been considered
in \cite{jeyak-cccq}.

\begin{theorem}\label{strog-dual-class-l} Let $S \subseteq X$ be a nonempty convex set, $f : X \rightarrow \overline \R$ a proper and convex function and $g : X \rightarrow Z^{\bullet}$ a
proper and $C$-convex function. If one of the regularity conditions $(RC^L_i)$, $i\in\{1,2,3,4,5\}$, is fulfilled, then $v(P_L)=v(D_L)$ and $(D_L)$ has an optimal solution.\end{theorem}

\begin{remark}\label{impl-class-l}\rm In case $X$ and $Z$ are Fr\' echet spaces, $S$ is a nonempty convex and closed set, $f : X \rightarrow \overline \R$ is a proper, convex and lower semicontinuous function
and $g : X \rightarrow Z^{\bullet}$ is a proper, $C$-convex and
$C$-epi closed function we have the following relations between
the above regularity conditions (cf. Remark \ref{impl-class})
$$(RC_1^L)\Rightarrow(RC_2^L)\Leftrightarrow(RC_3^L)\Rightarrow(RC_4^L)\Leftrightarrow(RC_5^L).$$
\end{remark}

We suppose in the following that $v(P_L)\in\R$. By some algebraic manipulations we get
$${\cal{E}}_{v(P_L)}:=\pr\nolimits_{Z \times\R}\big(\epi(\Phi_L-v(P_L))\big)$$
$$=\{(g(x)+z, f(x)-v(P_{L}) + \varepsilon):x\in \dom f \cap S \cap \dom g, z\in
C, \varepsilon \geq 0\}.$$
The set $-{\cal{E}}_{v(P_L)}$ is in analogy to the \emph{conic extension}, a notion used by Giannessi in the theory of \emph{image space analysis} (see \cite{giannessi}).

By means of the general scheme we can introduce now regularity conditions expressed by means of the quasi interior and quasi-relative interior for the primal-dual pair $(P_L)-(D_L)$, too. Besides the conditions $(RC_6^{\Phi_L})$ and $(RC_7^{\Phi_L})$, which in this case are exactly $(RC^L_6)$ and, respectively, $(RC^L_7)$ below, we consider like in the previous section a further one, which we denote by $(RC^L_{6'})$:
\begin{center}
\begin{tabular}{r|l}
$(RC^F_{6'})$ \ & \ $\exists x'\in\dom f \cap S$ such that $g(x')\in -\qri(C), \cl(C-C)=Z$ and\\
& \ $(0,0)\notin\qri\Big[\co\Big({\cal{E}}_{v(P_L)} \cup\{(0,0)\}\Big)\Big]$;
\end{tabular}
\end{center}

\begin{center}
\begin{tabular}{r|l}
$(RC^L_6)$ \ & \ $0\in\qi\big(g(\dom f \cap S \cap \dom g) +
C\big)$ and $(0,0)\notin\qri\Big[\co\Big({\cal{E}}_{v(P_L)}
\cup\{(0,0)\}\Big)\Big]$
\end{tabular}
\end{center}

and

\begin{center}
\begin{tabular}{r|l}
$(RC^L_7)$ \ & \ $0\in\qi\big[\big(g(\dom f \cap S
\cap \dom g) + C\big)-\big(g(\dom f \cap S \cap \dom g) + C\big)\big]$,  \\
& \ $0\in\qri\big(g(\dom f \cap S \cap \dom g) + C\big)$ and
$(0,0)\notin\qri\Big[\co\Big({\cal{E}}_{v(P_L)}
\cup\{(0,0)\}\Big)\Big]$.
\end{tabular}
\end{center}

Before studying the relations between these regularity conditions
we would like to notice that the set ${\cal{E}}_{v(P_L)}$ is
convex if and only if the function $(f,g): X \rightarrow \B \times
Z^\bullet$ is \emph{convexlike} on $\dom f \cap S \cap \dom g$
with respect to the cone $\R_+ \times C$, that is the set
$(f,g)(\dom f \cap S \cap \dom g) + \R_+ \times C$ is convex. This
property also implies that both sets $f(\dom f \cap S \cap \dom g)
+ \R_+$ and $g(\dom f \cap S \cap \dom g) + C$ are convex, too,
while the reverse implication does not always hold. Obviously, if
$S$ is a convex set, $f$ is a convex function and $g$ is a
$C$-convex function, then $(f,g)$ is convexlike on $\dom f \cap S
\cap \dom g$ with respect to the cone $\R_+ \times C$.

\begin{proposition}\label{comp-qri-l}{\it Let  $S \subseteq X$ be a nonempty set and $f : X \rightarrow \overline \R$ and $g : X \rightarrow
Z^{\bullet}$ be proper functions such that $v(P_L) \in \R$ and
$(f,g): X \rightarrow \B \times Z^\bullet$ is convexlike on $\dom
f \cap S \cap \dom g$ with respect to the cone $\R_+ \times C$
(the latter is the case if for instance $S$ is a convex set, $f$
is a convex function and $g$ is a $C$-convex function). The
following statements are true:
\begin{enumerate}\item[(i)] $(RC^L_{6'})\Rightarrow(RC^L_6)\Leftrightarrow(RC^L_7)$;

\item[(ii)] if $(P_L)$ has an optimal solution, then
$(0,0)\notin\qri\Big[\co\Big({\cal{E}}_{v(P_L)}
\cup\{(0,0)\}\Big)\Big]$ can be equivalently written as
$(0,0)\notin\qri\big({\cal{E}}_{v(P_L)}\big)$;

\item[(iii)] if $0\in\qi\big[\big(g(\dom f \cap S \cap \dom g) +
C\big)-\big(g(\dom f \cap S \cap \dom g) + C\big)\big]$, then
$(0,0)\notin\qri\Big[\co\Big({\cal{E}}_{v(P_L)}
\cup\{(0,0)\}\Big)\Big]$ is equivalent to
$(0,0)\notin\qi\Big[\co\Big({\cal{E}}_{v(P_L)}
\cup\{(0,0)\}\Big)\Big]$.\end{enumerate}}
\end{proposition}

\noindent {\bf Proof.} In view of Proposition \ref{comp-qri-phi}
we have to prove only the implication
$(RC^L_{6'})\Rightarrow(RC^L_6)$. Let us suppose the condition
$(RC^L_{6'})$ is fulfilled. We apply Lemma \ref{alte-prop-qri}(i)
with $U:=C$ and $V:=-g(\dom f \cap S \cap \dom g) - C$, which are
both convex sets. As $C$ is a convex cone both assumptions in this
statement are fulfilled and, consequently, $0 \in \qi\big(g(\dom f
\cap S \cap \dom g) + C+C\big) = \qi\big(g(\dom f \cap S \cap \dom
g) + C\big).$ This means that $(RC^L_6)$ holds.\hfill{$\Box$}

\begin{remark}\label{rel-magnanti}\rm For a special instance of $(P_L)$, in \cite{qri-siam} the regularity conditions $(RC^L_i), i \in \{6',6,7\}$, have been deduced  from $(RC^F_i), i \in \{6',6,7\}$, respectively,
by employing an approach due to Magnanti (cf. \cite{Magnanti})
which provides a link between the Fenchel and Lagrange duality
concepts. Let us mention that in the same spacial instance of
$(P_L)$, the condition $(RC_{6'}^L)$ has been considered in
\cite{qri-jota} in the framework of real normed spaces.
\end{remark}

\begin{remark}\label{nice-looking-cond-l}\rm (a) The condition $0\in\qi\big(g(\dom f \cap S \cap \dom g) + C\big)$
implies relation $0\in\qi\big[\big(g(\dom f \cap S \cap \dom g) +
C\big)-\big(g(\dom f \cap S \cap \dom g) + C\big)\big]$ in
Proposition \ref{comp-qri-l}(iii).

(b) The following implication holds
$$(0,0)\in\qi\Big[\co\Big({\cal{E}}_{v(P_L)}
\cup\{(0,0)\}\Big)\Big] \Rightarrow 0\in\qi\big(g(\dom f \cap S
\cap \dom g) + C\big)$$ and in this case similar comments as in
Remark \ref{nice-looking-cond-phi}(b) can be made.\end{remark}

A direct consequence of Theorem \ref{str-dual-phi-qri} and
Proposition \ref{comp-qri-l}(i) is the following strong duality
result concerning the primal-dual pair $(P_L)-(D_L)$ (see also
\cite{qri-siam}).

\begin{theorem}\label{str-dual-qri-l} Let  $S \subseteq X$ be a nonempty set and $f : X \rightarrow \overline \R$
and $g : X \rightarrow Z^{\bullet}$ be proper functions such that
$v(P_L) \in \R$ and $(f,g): X \rightarrow \B \times Z^\bullet$ is
convexlike on $\dom f \cap S \cap \dom g$ with respect to the cone
$\R_+ \times C$ (the latter is the case if for instance $S$ is a
convex set, $f$ is a convex function and $g$ is a $C$-convex
function). Suppose that one of the regularity conditions
$(RC^L_i)$, $i\in\{6',6,7\}$, holds. Then $v(P_L)=v(D_L)$ and
$(D_L)$ has an optimal solution.\end{theorem}

The following example considered by Daniele and Giuffr\` e in \cite{Daniele-Giuffre} shows that if we renounce to the condition $(0,0)\not\in\qri\big[\co({\cal{E}}_{v(P_L)}\cup\{(0,0)\})\big]$, then the strong duality result may fail.

\begin{example}\label{ren1}\rm Let be $X=Z=\ell^2$, $C=\ell^2_+$ and $S=\ell^2$. Take $f: \ell^2 \rightarrow \R$,
$f(x)=\langle c,x\rangle$, where $c=(c_n)_{n\in\N}, c_n=1/n$ for
all $n\in\N$ and $g:\ell^2\rightarrow \ell^2,g(x)=-Ax$, where
$(Ax)_n=(1/2^n)x_n$ for all $n\in\N$. Then $T =\{x\in \ell^2:Ax\in
\ell^2_+\}=\ell^2_+$. It holds $\cl(\ell^2_+ - \ell^2_+) = \ell^2$
and $\qri (\ell^2_+)=\{x=(x_n)_{n\in\N}\in \ell^2:x_n>0 \ \forall
n\in\N\}\neq\emptyset$ and one can easily find an element
$\overline{x}\in \ell^2$ with $g(\overline{x})\in-\qri
(\ell^2_+)$. We also have that
$$v(P_L)=\inf_{x\in T}\langle c,x\rangle=0$$ and $x=0$ is an
optimal solution of the primal problem. On the other hand, for
$z^*=(z^*_n)_{n\in\N} \in C^*=\ell_+^2$ it holds
$$\inf\limits_{x\in S}\{f(x)+\langle
z^*,g(x)\rangle\}=\inf\limits_{x\in \ell^2}\{\langle
c,x\rangle+\langle z^*,g(x)\rangle\}$$
$$=\inf\limits_{x=(x_n)_{n\in\N}\in \ell^2}
\left(\sum\limits_{n=1}^{\infty}\frac{1}{n}x_n-\sum\limits_{n=1}^{\infty}z^*_n\frac{1}{2^n}x_n\right)
=\inf\limits_{(x_n)_{n\in\N}\in
\ell^2}\sum\limits_{n=1}^{\infty}\left(\frac{1}{n}-\frac{z^*_n}{2^n}\right)x_n$$
$$=\left\{
\begin{array}{ll}
0, & \mbox {if } z^*_n=\frac{2^n}{n} \ \forall n\in\N,\\
-\infty, & \mbox{otherwise}.
\end{array}\right.$$ Since $(2^n/n)_{n\in\N}$ does not
belong to $\ell^2$, we obtain $v(D_L)=-\infty$, hence strong duality fails.

Moreover, it is not surprising that strong duality does not holds,
since not all the conditions in $(RC_i^L)$, $i\in\{6',6,7\}$, are
fulfilled. This is what we show in the following, namely that
$(0,0)\in\qi\big({\cal{E}}_{v(P_L)}\big)$. Take an arbitrary
element $(x^*,r^*)\in N_{{\cal{E}}_{v(P_L)}}(0,0)$, with
$x^*=(x_n^*)_{n\in\N}\in\ell^2$ and $r^*\in\R$. Then we have
\begin{equation}\label{tetalambda-ex} r^*(\langle
c,x\rangle +\varepsilon)+\langle x^*,g(x)+z\rangle\leq 0 \ \forall
x\in \ell^2 \ \forall \varepsilon\geq 0 \ \forall z\in
\ell_+^2,\end{equation} that is
$$r^*\left(\sum\limits_{n=1}^{\infty}\frac{1}{n}x_n+\varepsilon\right)+\sum\limits_{n=1}
^{\infty} x^*_n\left(-\frac{1}{2^n}x_n +z_n\right)\leq 0$$
$$\forall x=(x_n)_{n\in \N}\in \ell^2 \ \forall \varepsilon\geq 0 \ \forall
z=(z_n)_{n\in\N}\in \ell^2_+.$$ Taking $\varepsilon=0$ and $z_n=0$
for all $n\in\N$ in the relation above we get
$$\sum\limits_{n=1}^{\infty}\left(r^*\frac{1}{n}-\frac{1}{2^n}
x^*_n\right)x_n\leq 0 \ \forall x=(x_n)_{n\in \N}\in \ell^2,$$
which implies $x^*_n=r^* (2^n/n)$ for all $n\in\N$. Since $x^*\in
\ell^2$, we must have $r^*=0$ and hence $x^*=0$. Thus
$N_{{\cal{E}}_{v(P_L)}}(0,0)=\{(0,0)\}$ and so
$(0,0)\in\qi\big({\cal{E}}_{v(P_L)}\big)$ (cf. Proposition
\ref{charact-qi}).\end{example}

The following example underlines the applicability of the new
regularity conditions in opposition to the classical generalized
interior-point ones.

\begin{example}\label{ex-qri-l} \rm Consider again $X=Z=\ell^2(\N)$ and $C=\ell^2_+$. Take further $S=\ell^2_+$,
$f:\ell^2\rightarrow \R,$ $f(x)=\langle c,x\rangle$ and
$g:\ell^2\rightarrow\ell^2$, $g(x)=x-x^0$, where
$c,x^0\in\ell^2_+$ are arbitrary chosen such that $x^0_n>0$ for
all $n\in\N$. The feasible set of the primal problem is
$T=\ell^2_+\cap(x^0-\ell^2_+)\neq\emptyset$ and it holds
$$v(P_L)=\inf_{x\in T}\langle c,x\rangle=0,$$ while
$\overline{x}=0$ is an optimal solution of $(P_L)$. The condition
$\cl(C-C)=\ell^2$ is obviously satisfied and we have that (cf.
Example \ref{qri-lp+}) $\{x\in \dom f \cap S:
g(x)\in-\qri(C)\}=\{x=(x_n)_{n\in\N}\in\ell^2:0\leq x_n < x^0_n \
\forall n\in\N\}$. This is a nonempty set, hence the Slater type
condition is also fulfilled. We prove in the following that
$(0,0)\not\in\qi\big({\cal{E}}_{v(P_L)}\big)$. An arbitrary
element $(x^*,r^*)\in \ell^2\times\R$ belongs to
$N_{{\cal{E}}_{v(P_L)}}(0,0)$ if and only if $$r^*(\langle
c,x\rangle+\varepsilon)+\langle x^*,x-x^0+z\rangle\leq 0 \ \forall
x\in\ell^2_+ \ \forall\varepsilon\geq 0 \ \forall z\in\ell^2_+.$$
One can observe that $(0,-1)\in N_{{\cal{E}}_{v(P_L)}}(0,0)$,
which ensures that $N_{{\cal{E}}_{v(P_L)}}(0,0)\neq\{(0,0)\}$. By
using Proposition \ref{charact-qi} we obtain the conclusion. An
application of Proposition \ref{comp-qri-l}(iii) yields that
$(RC^L_{6'})$ is fulfilled (and, consequently, also $(RC^L_i)$, $i
\in \{6,7\}$) and hence strong duality holds (cf. Theorem
\ref{str-dual-qri-l}). On the other hand, since $g(\dom f \cap S
\cap \dom g)+C=\ell^2_+-x^0$, none of the regularity conditions
$(RC^L_i)$, $i \in \{1, 2, 3, 4, 5\}$, can be applied to this
problem (see Example \ref{qri-lp+}). Notice also that the optimal
objective value of the Lagrange dual problem is
$$v(D_L)=\sup_{z^*\in\ell^2_+}\inf_{x\in\ell^2_+}\{\langle
c,x\rangle+\langle z^*,
x-x^0\rangle\}$$$$=\sup_{z^*\in\ell^2_+}\Big\{-\langle
z^*,x^0\rangle+\inf_{x\in\ell^2_+}\langle
c+z^*,x\rangle\Big\}=\sup_{z^* \in\ell^2_+}\langle
-z^*,x^0\rangle=0$$ and $\bar z^*=0$ is an optimal solution of the
dual.\end{example}

Next we show that in general the condition $(RC_6^L)$ (and
implicitly also $(RC_7^L)$) is weaker than $(RC_{6'}^L)$.

\begin{example}\label{ex-rc7-rc6-l} \rm We work in the following setting: $X=Z=\ell^2(\R)$, $C=\ell^2_+(\R)$,
$S=\ell^2_+(\R)$, while the functions $f:\ell^2(\R)\rightarrow
\R$, $g:\ell^2(\R)\rightarrow\ell^2(\R)$ are defined as
$f(s)=\|s\|$ and $g(s)=-s$ for all $s\in\ell^2(\R)$, respectively.
For the primal problem we have
$$v(P_L)=\inf_{x\in\ell^2_+(\R)}\|s\|=0$$ and $s=0$ is the unique optimal
solution of $(P_L)$. Since $\qri\big(\ell^2_+(\R)\big)=\emptyset$
(cf. Example \ref{qri-vid}), the condition $(RC_{6'}^L)$ fails.
Further, $g(\dom f \cap S \cap \dom
g)+C=-\ell^2_+(\R)+\ell^2_+(\R)=\ell^2(\R)$, hence
$0\in\qi\big(g(\dom f \cap S \cap \dom g)+C\big)$. Now one can
easily prove like in Example \ref{ex-qri-l} that
$(0,0)\not\in\qi\big({\cal{E}}_{v(P_L)}\big)$, thus the condition
$(RC_6^L)$ is fulfilled (cf. Proposition \ref{comp-qri-l}(iii)),
hence strong duality holds. The optimal objective value of the
dual problem is
$$v(D_L)=\sup_{z^*
\in\ell^2_+(\R)}\inf_{s\in\ell^2_+(\R)}\{\|s\|-\langle
z^*,s\rangle\}.$$ For every $z^* \in\ell^2_+(\R)$ the inner
infimum in the above relation can be written as (cf. \cite[Theorem
2.8.7]{Zal-carte}) $$\inf_{s\in\ell^2_+(\R)}\{\|s\|-\langle
z^*,s\rangle\}=-\sup_{s\in\ell^2_+(\R)}\{\langle
z^*,s\rangle-\|s\|\}=
-(\|\cdot\|+\delta_{\ell^2_+(\R)})^*(z^*)$$$$=-(\delta_{\overline{B}(0,1)}\Box
\delta_{-\ell^2_+(\R)})(z^*)=-\delta_{\overline{B}(0,1)-\ell^2_+(\R)}(z^*),$$
where $\overline{B}(0,1)$ is the closed unit ball of
$\big(\ell^2(\R)\big)^*=\ell^2(\R)$. We get $v(D_L)=0$ and every
$z^*\in\ell^2_+(\R)\cap\big(\overline{B}(0,1)-\ell^2_+(\R)\big)$
is an optimal solution of the dual (in particular also $\bar
z^*=0$).\end{example}

From Proposition \ref{strong-impl-qi-phi} we obtain the next result.

\begin{proposition}\label{strong-impl-qi-l} Suppose that for the primal-dual pair $(P_L)-(D_L)$ strong duality holds.
Then $(0,0)\notin\qri\Big[\co\Big({\cal{E}}_{v(P_L)} \cup\{(0,0)\}\Big)\Big]$. \end{proposition}

This results allows giving the following extended  scheme involving the regularity conditions for the primal-dual pair $(P_L)-(D_L)$.

\begin{proposition}\label{comp-qri-class-l} Suppose that $X$ and $Z$ are  Fr\'{e}chet spaces, $S$ is a nonempty
convex and closed set, $f : X \rightarrow \overline \R$ is a
proper, convex and lower semicontinuous function and $g : X
\rightarrow Z^{\bullet}$ is a proper, $C$-convex and $C$-epi
closed function. The following relations hold
$$(RC_1^L)\Rightarrow(RC_2^L)\Leftrightarrow(RC_3^L)\Rightarrow(RC_6^F)\Leftrightarrow(RC_7^F).$$\end{proposition}

\begin{remark}\label{rc1-impl-rc6-l}\rm One can notice that the implications $$(RC_1^L)\Rightarrow(RC_6^L)\Leftrightarrow(RC_7^L)$$ hold in the framework of separated locally convex spaces and for
$S$ a convex set, $f : X \rightarrow \overline \R$ a proper and
convex function and $g : X \rightarrow Z^{\bullet}$ a proper and
$C$-convex function (nor completeness for the spaces involved
neither topological assumptions for the functions considered are
needed here).\end{remark}

In general the conditions $(RC_i^L)$, $i\in\{4,5\}$, cannot be compared with $(RC_i^L)$, $i\in\{6',6,7\}$. Example \ref{ex-qri-l} provides a situation when $(RC_i^L)$, $i\in\{6',6,7\}$, are fulfilled, unlike $(RC_i^L)$, $i\in\{4,5\}$. In the following example the conditions $(RC_i^L)$, $i\in\{4,5\}$, are fulfilled, while $(RC_i^L)$, $i\in\{6',6,7\}$, fail.

\begin{example}\label{sqri-qri-nu-l}\rm We use the idea from Example \ref{sqri-qri-nu} and consider
$(X,\|\cdot\|)$ a nonzero real Banach space and $x_0^*\in
X^*\setminus\{0\}$. Let further be  $Z= X$, $C=\{0\}$, $S= \Ker
x_0^* \subseteq X$, $f : X \rightarrow \R$, $f=\|\cdot\|$  and
$g:X \rightarrow X$, $g(x)=x$ for all $x \in X$. The set $S$ is
convex and closed, $f$ is convex and continuous, while $g$ is
obviously $C$-convex and $C$-epi closed. Moreover, $v(P_L) = 0$
and $\bar x=0$ is the unique optimal solution of $(P_L)$. Further,
$g(\dom f \cap S \cap \dom g) + C =\Ker x_0^*$, which is a closed
linear subspace of $X$, hence $(RC_i^L)$, $i\in\{4,5\}$, are
fulfilled. Further, $\big(g(\dom f \cap S \cap \dom g) +
C\big)-\big(g(\dom f \cap S \cap \dom g) + C\big) =g(\dom f \cap S
\cap \dom g) + C=\Ker x_0^*$ and, as we have seen that one cannot
have $\Ker x_0^*=X$, all the three regularity conditions
$(RC_i^F)$, $i\in\{6',6,7\}$, fail.

Notice that $v(D_L) = \sup_{z^* \in X^*} -(\|\cdot\| + \delta_{\Ker x_0^*})^*(-z^*) = \sup_{z^* \in X^*} -\delta_{B_*(0,1)+\R x_0^*}(-z^*) = 0$ and that the set of optimal solutions of $(D_L)$
coincides with $B_*(0,1)+\R x_0^*$.
\end{example}

\begin{remark}\label{closedness-type-l}\rm
It is worth mentioning that for the primal-dual pair $(P_L)-(D_L)$
besides the generalized interior-point regularity conditions also
\emph{closedness-type regularity conditions} have been considered
in the literature. In the general setting considered in this
section we have the following condition of this type (cf.
\cite{b-hab})
\begin{center}
\begin{tabular}{r|l}
$(RC^L_8)$ \ & \ $S$ is closed, $f$ is lower semicontinuous, $g$ is $C$-epi closed and\\
& \ $\cp\limits_{z^*\in C^*} \epi(f+(z^*g) + \delta_S)^*$ is
closed in $(X^*,w(X^*,X))\times\R$.
\end{tabular}
\end{center}
Whenever $S$ is a nonempty convex set, $f$ is a proper and convex
function and $g$ is a proper and $C$-convex function with $\dom f
\cap S \cap g^{-1}(-C) \neq \emptyset$ and $(RC^L_8)$ is fulfilled
it holds
\begin{equation}\label{stab-str-l}(f + \delta_T)^*(x^*)=\min\{(f+(z^*g) + \delta_S)^*(x^*) :z^*\in C^*\} \ \forall x^*\in X^*\end{equation}
and this obviously guarantees strong duality for $(P_L)-(D_L)$.
When, additionally, $S$ is closed, $f : X \rightarrow \overline
\R$  is lower semicontinuous  and $g : X \rightarrow Z^{\bullet}$
is $C$-epi closed, then $(RC^L_8)$ is fulfilled if and only if
\eqref{stab-str-l} holds (cf. \cite[Theorem 1]{BGW}).

In case $X$ and $Z$ are  Fr\'{e}chet spaces, $S$ is a nonempty convex and closed set, $f : X \rightarrow \overline \R$ a proper, convex and lower semicontinuous function and $g : X \rightarrow Z^{\bullet}$ a proper, $C$-convex and $C$-epi closed function we have the following relations between the regularity conditions considered for the primal-dual pair $(P_L)-(D_L)$ (cf. \cite{b-hab}) $$(RC_1^L)\Rightarrow(RC_2^L)\Leftrightarrow(RC_3^L)\Rightarrow(RC_4^L)\Leftrightarrow(RC_5^L)\Rightarrow(RC_8^L),$$ in general the implications above being strict. Notice that for the implication $(RC_1^F)\Rightarrow(RC_8^F)$ one can omit assuming completeness for $X$ and $Z$.

Finally, we want to point out the fact that in general $(RC_6^L)$ (and automatically also $(RC_7^L)$) and $(RC_8^L)$ are not comparable. This assertion can be illustrated by adapting the examples
\ref{cl-nu-qri-f-2} and \ref{nu-cl-qri-da}.
\end{remark}

\begin{remark}\label{alte-cond-qri-lit}\rm One can find in the literature various results, the majority of them recently published, 
where regularity conditions employing the quasi interior and the quasi-relative interior have been provided. Unfortunately, these statements
have either superfluous, or contradictory hypotheses. One can overcome these drawbacks by using instead the results
presented in this section.

A regularity condition for strong duality for the pair
$(P_L)-(D_L)$ was proposed by Cammaroto and Di Bella in
\cite[Theorem 2.2]{Cammaroto}. However, in \cite[Section
4]{qri-siam} it has been shown that this theorem finds no
application, since the hypotheses considered by the authors are contradictory.

Daniele, Giuffr\`{e}, Idone and Maugeri formulated in
\cite{Daniele-Giuffre-Maugeri} a strong duality theorem for the
primal-dual pair $(P_L)-(D_L)$. Besides some conditions involving
the quasi-relative interior, the authors considered a supplementary condition,
called \emph{Assumption S}. We emphasize that \emph{Assumption
S} is fulfilled if and only for the pair $(P_L)-(D_L)$ strong duality holds (see
\cite[Corollary 2.1]{qri-jota}). This means that there is no need to impose the other assumptions,
which are expressed by using the quasi-relative interior. In this way the aim of proposing sufficient
conditions for duality followed by the authors in the mentioned paper is not attained.\end{remark}

\end{document}